\documentclass[10pt]{article}
\usepackage{amsmath,amssymb}

\newcommand{\pf }{{\noindent{\bf Proof: }}}

\newcommand{\bra}{\langle}
\newcommand{\ket}{\rangle}
\newcommand{\ep}{\hfill {$\Box$}}

\newtheorem{thm}{Theorem}[section]
\newtheorem{cor}{Corollary}[section]
\newtheorem{defin}{Definition}[section]
\newtheorem{lem}{Lemma}[section]
\newtheorem{prop}{Proposition}[section]

\title{Commutation Relations for Unitary Operators II}
\author{\small{M.A. Astaburuaga, O. Bourget, V.H. Cort\'es \footnote{Supported by the Grant Fondecyt 1120786, ECOS-Conicyt C10E10}}\\
\small{Facultad de Matem\'aticas, Pontificia Universidad Cat\'olica de Chile,}\\
\small{Av. Vicu\~na Mackenna 4860, Macul, Santiago, Chile}\\
\small{E-mail: bourget@mat.puc.cl}\\
\small{phone: (56 2) 2354 4509}\\
}
\date{}
\begin{document}
\maketitle

\begin{abstract}
Let $f$ be a regular non-constant symbol defined on the $d$-dimensional torus ${\mathbb T}^d$ with values on the unit circle. Denote respectively by $\kappa$ and $L$, its set of critical points and the associated Laurent operator on $l^2({\mathbb Z}^d)$. Let $U$ be a suitable unitary local perturbation of $L$. We show that the operator $U$ has finite point spectrum and no singular continuous component away from the set $f(\kappa)$. We apply these results and provide a new approach to analyze the spectral properties of GGT matrices with asymptotically constant Verblunsky coefficients. The proofs are based on positive commutator techniques. We also obtain some propagation estimates.
\end{abstract}

\noindent{\it Keywords:} Spectrum, Commutator, Unitary Laurent operator, GGT matrices, perturbation.

\section{Introduction}

The spectral analysis of unitary operators on Hilbert spaces is naturally related to the study of dynamical systems \cite{N2}, periodic time-dependent quantum dynamical systems \cite{ev} and the theory of orthogonal polynomials on the unit circle e.g. \cite{S1}, \cite{S2}. That analysis can be achieved by means of commutation relations techniques. We refer to the introduction of \cite{abc1} for a bibliographical discussion.

In this paper, we apply such techniques to analyze the effects of local and regular perturbations on the spectral properties of some unitary Laurent operators defined on  $l^2({\mathbb Z}^d)$ (Theorems \ref{laurent-glob} and \ref{laurent-loc}). These results allows to recover the spectral properties of GGT matrices with asymptotically constant Verblunsky coefficients (Theorem \ref{ggtperturbed1}). This is the first time these techniques are applied within the field of orthogonal polynomials on the unit circle. An extension of this approach to a wider class of Verblunsky coefficients in the spirit of \cite{sah1}, \cite{sah2}, \cite{sah3} remains to be done. In order to provide a dynamical interpretation to those results, we have also explicited some related propagation estimates (Propositions \ref{laurent-2ub} and \ref{laurent-2lb})

The manuscript is structured as follows. The main spectral results are introduced in Section 2 and applied in Section 3. The proofs are developed in Section 4. The propagation estimates are considered in Section 5. Some complementary tools are gathered in Section 6.  \\

\noindent {\bf Notations:} ${\cal H}$ denotes a infinite-dimensional (complex) Hilbert space and ${\cal B}({\cal H})$ the algebra of bounded operators acting on ${\cal H}$. The resolvent set of an operator $B\in {\cal B}({\cal H})$ is denoted by $\rho(B)$ and its spectrum by $\sigma(B) := {\mathbb C}\setminus \rho(B)$. The unit circle, the open unit disc and the one-dimensional torus ${\mathbb R}/2\pi {\mathbb Z}$ are denoted by $\partial {\mathbb D}$, ${\mathbb D}$ and ${\mathbb T}$ respectively. To any bounded Borel function $\Phi$ on $\partial {\mathbb D}$ is associated a unique function $\phi$ defined on ${\mathbb T}$ by: $\phi(\theta)=\Phi(e^{i\theta})$, for all $\theta \in {\mathbb T}$. If $U$ is a unitary operator defined on ${\cal H}$ and if its spectral family is denoted by $(E_{\Delta})_{\Delta \in {\cal B}({\mathbb T})}$, where ${\cal B}({\mathbb T})$ stands for the family of Borel sets of ${\mathbb T}$, we will have that:
\begin{equation*}
\Phi(U)=\int_{\mathbb T}\phi(\theta)dE_{\theta}=\int_{\mathbb T}\Phi(e^{i\theta})dE_{\theta}\enspace .
\end{equation*}
The continuous and point subspaces of the operator $U$ are respectively denoted by ${\cal H}_c(U)$ and ${\cal H}_{pp}(U)$. In this paper, we deal mostly with the Hilbert spaces $L^{2} ({\mathbb T}^d)$ and $l^{2} ({\mathbb Z}^d)$, where $d\in {\mathbb N}$. $L^{2} ({\mathbb T}^d)$ is the Hilbert space of square integrable complex functions with inner product
\begin{equation*}
\langle f, g  \rangle =  \frac{1}{(2\pi)^d  } \,  \int_{{\mathbb T}^d} \overline{f(\theta)} g(\theta) \, d\theta \enspace ,
\end{equation*}
while $l^{2} ({\mathbb Z}^d)$ denotes the Hilbert space of square summable complex functions defined on the lattice $\mathbb Z^d$, $(\varphi_{\alpha})_{\alpha \in {\mathbb Z}}$, with inner product $\langle (\varphi_{\alpha})_{\alpha \in {\mathbb Z}^d }, (\psi_{\alpha})_{\alpha \in {\mathbb Z}^d } \rangle = \sum_{\alpha \in  {\mathbb Z}^d} \overline{\varphi_{\alpha}} \,  \psi_{\alpha}$. Given $\beta \in {\mathbb Z}^d$, we define the function $e_{\beta}$ on ${\mathbb Z}^d$ by $(e_{\beta})_{\alpha} = \delta_{\alpha \beta}$. $(e_{\beta})_{\beta \in {\mathbb Z}^d}$ is the canonical orthonormal basis of $l^{2} ({\mathbb Z}^d)$. With these notations, the Fourier transform  ${\mathcal F} : L^2({\mathbb T}^d) \to l^{2} ({\mathbb Z}^d)$ is defined by: ${\mathcal F} f = (\hat{f}_{\alpha} )_{\alpha \in {\mathbb Z}^d  }$ where $  \hat{f}_{\alpha} $ is the $\alpha$-Fourier coefficient of the function $f\in L^{2} ({\mathbb T}^d)$:
 \begin{equation}\label{fourier1}
 \hat{f}_{\alpha} = \frac{1}{(2 \pi)^d} \int_{{\mathbb T}^d} e^{-i\theta \cdot \alpha} f(\theta) \, d\theta \enspace .
\end{equation}
The shift operators $(T_j)_{j\in \{1,\ldots,d\}}$ are examples of unitary operators defined on $l^{2} ({\mathbb Z}^d)$. Their action on the orthonormal basis of $l^2({\mathbb Z}^d)$ is given by: for all $\beta \in {\mathbb Z}^d$, $T_j e_{\beta} = e_{S_j(\beta)}$, where $S_j(\beta) = (\beta_1, \ldots, \beta_{j-1}, \beta_j+1, \beta_{j+1}, \ldots , \beta_d)$. For all $(i,j)\in \{1,\ldots, d\}^2$: $T_j^*=T_j^{-1}$ and $[T_i,T_j]=0$. For any $\alpha=(\alpha_1,\ldots,\alpha_d)\in {\mathbb Z}^d$, we write $T^{\alpha}=T_1^{\alpha_1} \ldots T_d^{\alpha_d}$.

The spaces $L^{\infty} ({\mathbb T}^d)$, $C^0({\mathbb T}^d)$ and $C^k({\mathbb T}^d)$ ($k\in {\mathbb N}$) stand respectively for the linear spaces of essentially bounded complex functions, continuous complex functions and $k$-th continuously differentiable complex functions defined on ${\mathbb T}^d$. The Wiener algebra is defined by: ${\cal A}({\mathbb T}^d):= \{f\in L^{\infty}({\mathbb T}^d); (\hat{f}_{\alpha})_{\alpha \in {\mathbb Z}^d}\in l^1({\mathbb Z}^d)\}$ \cite{ed}, \cite{katz}.

\section{Preliminaries and Abstract Results}

Let $d\in {\mathbb N}$ and denote by $(X_j)_{j=1}^d$ the family of linear operators defined on the canonical orthonormal basis of $l^2({\mathbb Z}^d)$ by: $X_j e_{\beta} = \beta_j e_{\beta}$. The operators $(X_j)_{j=1}^d$ are essentially self-adjoint on the linear span of $(e_{\beta})_{\beta \in {\mathbb Z}^d}$. Their self-adjoint extensions are also denoted by $X_j$. For all $j\in \{1,\ldots, d\}$, ${\cal F}^*X_j {\cal F}= -i\partial_{\theta_j}$. We define the set
\begin{equation*}
{\cal D}_X=\{\varphi \in l^2({\mathbb Z^d}); \sum_{\alpha \in {\mathbb Z}^d}(1+\alpha_1^2+\ldots + \alpha_d^2) |\varphi_{\alpha}|^2 < \infty\}\, .
\end{equation*}
Note that ${\cal D}_X=\cap_{j=1}^d {\cal D}(X_j)={\cal D}(\bra X\ket)$ where $\bra X \ket:= \sqrt{X_1^2+\ldots+X_d^2 +1}$. In Sections 2.1, 2.2 and 2.3, we introduce the ingredients which allows to state Theorems \ref{laurent-glob}, \ref{laurent-loc} and Proposition \ref{transl}.

\subsection{Laurent Operators}

Let $f$ be a function in $L^{\infty}({\mathbb T}^d)$ and denote the sequence of its Fourier coefficients by $(\hat{f}_{\alpha})_{\alpha \in {\mathbb Z}^d}$. The Laurent operator $L_f$ associated to the symbol $f$ is the bounded discrete convolution operator defined for all $\psi\in l^2({\mathbb Z}^d)$  by
\begin{equation*}
L_f  \psi :=  \mathcal{F}(f)*\psi \enspace,
\end{equation*}
that is for any $\beta \in {\mathbb Z}^d,$
$$
(L_f \psi)_{\beta} =   \sum_{\alpha \in {\mathbb Z}^d} \hat{f}_{\alpha} \psi_{\beta-\alpha} =\sum_{\alpha \in {\mathbb Z}^d} \hat{f}_{\alpha}  (T^{\alpha}\psi)_{\beta}\enspace .
$$
The bounded linear operator $L_f$ is unitarily equivalent to the operator multiplication by $f$ on $L^2({\mathbb T}^d)$, denoted by $M_f$ or $f(\cdot)$: $L_f ={\cal F} f(\cdot) {\cal F}^*$. If the symbol $f$ belongs to the Wiener algebra ${\cal A}({\mathbb T}^d)$, then $L_f$ rewrites as the following norm convergent series:
\begin{equation*}
L_f=\sum_{\alpha \in {\mathbb Z}^d} \hat{f}_{\alpha} T^{\alpha} \enspace .
\end{equation*}
For any functions $f$ and $g$ in $L^{\infty}({\mathbb T}^d)$ and any $c\in {\mathbb C}$, we have that: $L_{f+g}=L_f+L_g$, $L_{fg}=L_f L_g$, $L_{cf}=cL_f$, $L_f^*=L_{\bar{f}}$. In particular, $[L_f,L_g]=0$ and $L_1=I$. For a real-valued function $f$ in $L^{\infty}({\mathbb T}^d)$, the Laurent operator $L_f$ is self-adjoint. If $f\in L^{\infty}({\mathbb T}^d)$ has values on $ \partial {\mathbb D}$ ($|f|=1$), $L_f$ is unitary.

\noindent{\bf Remark:} If $f$ is a continuously differentiable function with values in $\partial {\mathbb D}$, then for all $j\in \{1,\ldots,d\}$, $\Re(f\partial_{\theta_j}\bar{f})=0$. In other words, the Laurent operators $(L_{if\partial_{\theta_j} \bar{f}})_{j\in \{1,\ldots,d\}}$ are self-adjoint.

It is also well-known that the spectral properties of $L_f$ are related to the properties of the function $f$ (see e.g. paragraph 7.1.4 in \cite{abmg}). If $f\in  L^{\infty}({\mathbb T}^d)$, $|f|=1$, we have that:
\begin{itemize}
\item the spectrum of $L_f$ coincides with the essential range of $f$. If in addition $f$ is continuous on ${\mathbb T}^d$, we have that $\sigma(L_f)=f({\mathbb T}^d)=$ Ran $f$, which is a connected and compact subset of $\partial {\mathbb D}$.
\item $\lambda$ is an eigenvalue of $L_f$ if and only if $f^{-1}(\{\lambda\})$ has non zero Lebesgue measure.
\item $L_f$ has purely absolutely continuous spectrum in a subset $e^{i\Theta}$ of $\partial {\mathbb  D}$ ($\Theta\subset {\mathbb T}$) if and only if for any Borel set $N\subset e^{i\Theta}$ of zero Lebesgue measure, $f^{-1}(N)$ is also of measure zero.
\item $L_f$ has non-trivial singular continuous spectrum if and only if there exists a Borel set $N\subset \partial {\mathbb D}$ of zero Lebesgue measure, such that $f^{-1}(N)$ has non-zero measure but $f^{-1}(\{\lambda\})$ is of zero measure for each $\lambda\in N$.
\end{itemize}

Lastly, we denote the set of critical points of the symbol $f$ by: $\kappa_f = \{\theta \in {\mathbb T}^d; f$ is not differentiable at $\theta$ or $\nabla f(\theta)=0\}$. If $f$ belongs to $C^1({\mathbb T}^d)$ with values on $\partial {\mathbb D}$, the sets $\kappa_f$ and $f(\kappa_f)$ are compact subsets of ${\mathbb T}^d$ and $\sigma (L_f)=\mathrm{Ran} f \subset \partial {\mathbb D}$ respectively.

\subsection{Commutators and Regularity}

Let ${\cal H}$ be a Hilbert space, $B\in {\cal B}({\cal H})$ and $A$ be a self-adjoint operator (densely) defined on ${\cal H}$ with domain ${\cal D}(A)$. The operator $B$ is of class $C^1$ with respect to $A$ if the sesquilinear form defined on ${\cal D}(A)\times {\cal D}(A)$ by $(\varphi,\psi)\mapsto \langle A\varphi,B\psi\rangle - \langle \varphi,BA\psi\rangle$ is continuous for the topology induced by ${\cal H}\times {\cal H}$ i.e: there exists $C>0$ such that for all $(\varphi,\psi)\in {\cal D}(A)\times {\cal D}(A),$
\begin{equation}\label{qform}
\left| \langle A\varphi,B\psi\rangle - \langle \varphi,BA\psi\rangle \right| \leq C \|\varphi\| \|\psi\| \enspace .
\end{equation}
The form extends continuously as a bounded sesquilinear form on ${\cal H}\times {\cal H}$. The bounded linear operator associated to that extension is denoted by $\mathrm{ad}_A (B) =[A,B]$. We also write: $C^1(A):=\{B\in {\cal B}({\cal H}); B$ is of class $C^1$ w.r.t $A\}$.

\noindent{\bf Remark:} For practical purposes, it is enough to check estimate (\ref{qform}) on ${\cal S}\times {\cal S}$ where ${\cal S}$ is a core for $A$. For alternative characterizations of the class $C^1(A)$, see Section 6.2.

With the convention that $C^0(A):={\cal B}({\cal H})$ and $\mathrm{ad}^0_A B:=B$, higher order commutation relations are characterized inductively for $n\in {\mathbb N}$ as follows: $B\in {\cal B}({\cal H})$ is of class $C^n$ with respect to $A$ if $B\in C^{n-1}(A)$ and $\mathrm{ad}^{n-1}_A B \in C^1(A)$. Accordingly, we write $\mathrm{ad}_A^n B := \mathrm{ad}_A (\mathrm{ad}_A^{n-1} B)$, $C^n(A):=\{B\in {\cal B}({\cal H}); B$ is of class $C^n$ w.r.t $A\}$ and $C^{\infty}(A):= \cap_{n\in {\mathbb N}}C^n(A)$ (see Section 6.1 for more references).

Fractional regularity scales can also be considered (see \cite{abmg}, \cite{sah} and also Definition \ref{csp}). For the moment, let us just mention that for $B\in {\cal B}({\cal H})$ we say that:
\begin{itemize}
\item $B\in {\cal C}^{0,1}(A)$ if:
\begin{equation*}
\int_0^1 \|e^{iA\tau}B e^{-iA\tau}-B\|\,\frac{d\tau}{\tau} < \infty \enspace .
\end{equation*}
\item $B\in {\cal C}^{1,1}(A)$ if:
\begin{equation*}
\int_0^1 \|e^{iA\tau}B e^{-iA\tau}+e^{-iA\tau}B e^{iA\tau} -2B\|\,\frac{d\tau}{\tau^2} < \infty \enspace .
\end{equation*}
\end{itemize}
Clearly, ${\cal C}^{0,1}(A)$ and ${\cal C}^{1,1}(A)$ are linear subspaces of ${\cal B}({\cal H})$, stable under adjunction $^*$. It is also know that if $B\in C^1(A)$ and $\mathrm{ad}_A B\in {\cal C}^{0,1}(A)$, then $B\in {\cal C}^{1,1}(A)$ and that $C^2(A)\subset {\cal C}^{1,1}(A)\subset C^1(A)$ (see inclusions 5.2.19 in \cite{abmg}).

\noindent{\bf Remark:} The regularity hypotheses on the symbol $f$ introduced in Section 2.1 can be reinterpreted in terms of some regularity properties of $L_f$ w.r.t. some suitable commutation operations. For example, $(T_i)_{i=1}^d\subset \cap_{j=1}^d C^{\infty}(X_j)$: given any $(i,j)\in \{1,\ldots,d\}^2$ and any nonnegative integer $l$, $\mathrm{ad}_{X_j}^l(T_i)=\delta_{ij}T_i$ and $\mathrm{ad}_{X_j}^l(T_i^*)=(-1)^l \delta_{ij} T_i^{*}$. By Fourier transform (\ref{fourier1}), we deduce that for $h\in C^n({\mathbb T}^d)$, $L_h\in \cap_{j=1}^d C^n(X_j)$ and
\begin{equation}\label{xL-Lx}
\mathrm{ad}_{X_j} L_h = L_{-i\partial_{\theta_j}h}=-i L_{\partial_{\theta_j}h} \enspace .
\end{equation}

In this paper, we shall also relate the commutation conditions described above with some regularity properties of the resolvent of unitary operators. We say that a limiting absorption principle (LAP) holds for a unitary operator $U$ on some open subset $\Theta\subset {\mathbb T}$, w.r.t. a self-adjoint operator $A$, when:
\begin{itemize}
\item For any compact subset $K\subset \Theta$
\begin{equation*}
\sup_{ |z| \neq 1, \arg z \in K} \|\bra A \ket^{-1}(1-zU^*)^{-1}\bra A \ket^{-1} \| < \infty \enspace .
\end{equation*}
\item If $z$ tends to $e^{i\theta} \in e^{i\Theta}$ (non-tangentially), then $\bra A \ket^{-1}(1-zU^*)^{-1}\bra A \ket^{-1}$ converges in norm to a bounded operator denoted $R^+(\theta)$ (resp. $R^-(\theta) )$ if $|z| <1$ (resp. $|z| >1$). This convergence is uniform on any compact subset $K\subset \Theta$.
\item The operator-valued functions defined by $R^{\pm}$ are continuous on $\Theta$, with respect to the norm topology of ${\cal B}({\cal H})$.
\end{itemize}
\noindent{\bf Remark:} If $e^{i\Theta}\subset \rho (U)$, the properties described above are trivially satisfied. In this case, $R^+(\theta)=R^-(\theta)=(1-e^{i\theta}U^*)^{-1}$ for all $\theta\in \Theta$.

\subsection{Conjugate Operators}

Let $g=(g_j)_{j=1}^d \subset C^2({\mathbb T}^d)$ a family of real-valued functions and denote $L_g = (L_{g_j})_{j=1}^d$. Since $L_{g_j}\in \cap_{i=1}^d C^1(X_i)$, $L_{g_j}{\cal D}_X\subset {\cal D}_X$ for all $j\in \{1,\ldots,d\}$. To such a family $g$, we associate a symmetric operator $A_g$ defined on ${\cal D}_X$ by:
\begin{eqnarray}
A_g & :=& \frac{1}{2}\left( L_g \cdot X + X \cdot L_g \right)= \frac{1}{2} \sum_{j=1}^d L_{g_j} X_j + X_j L_{g_j}   \label{ag}\\
&=& \frac{1}{2} {\cal F} \left(g \cdot (-i\nabla) + (-i\nabla) \cdot g\right) {\cal F}^*\nonumber
\end{eqnarray}
where we have used that $\mathcal{F}^* L_{g_j} \mathcal{F} = g_j (\cdot)=M_{g_j} $ and $ \mathcal{F}^* X_j \mathcal{F} :=-i\partial _{\theta_j}$, $j$ in $\{1,\ldots,d\}$. Following the proof of Proposition 7.6.3 (a) in \cite{abmg}, we can show that the operator $i(g \cdot \nabla + \nabla \cdot g)$ is (defined and) essentially self-adjoint on $C^2({\mathbb T}^d)$. As ${\cal F} C^2({\mathbb T}^d) \subset {\cal D}_X$, we have obtained that:
\begin{lem}\label{sa} Let $g=(g_j)_{j=1}^d \subset C^2({\mathbb T}^d)$ a family of real-valued functions. Then, the operator $A_g$ defined on ${\cal D}_X$ by (\ref{ag}) is essentially self-adjoint. Moreover, $A_g\bra X\ket^{-1}$ and $A_g^2\bra X\ket^{-2}$ are bounded.
\end{lem}
The self-adjoint extension of $A_g$ defined in Lemma \ref{sa} is also denoted $A_g$.

\noindent{\bf Remark:} For the canonical projections $(p_j)_{j=1}^d$ defined on ${\mathbb T}^d$ by $p_j(\theta)=\theta_j$, we have that: $A_{p_j}=X_j$.

As a direct application of the  Fourier transform (\ref{fourier1}), we have also that:
\begin{lem}\label{H0-4} Let $f\in C^0({\mathbb T}^d)$ and $(n_1, n_2) \in {\mathbb N}^2$ with $n_2 \geq 2$. Consider $g=(g_j)_{j=1}^d \subset C^{n_2}({\mathbb T}^d)$ a family of real-valued functions. Suppose that there exists an open set $\Theta \subset {\mathbb T}^d$, which contains the support of $g$, on which $f$ is of class $C^{n_1}$. Then $L_f\in C^1(A_g)$ and $\mathrm{ad}_{A_g} L_f = L_{-ig\cdot \nabla f} $. It follows that $\mathrm{ad}_{A_g} L_f \in C^{\min(n_1-1,n_2)}(A_g)$. In particular, if $n_1\geq 3$, $L_f \in C^{n_1}(A_{if \nabla\bar{f}})$ and $\mathrm{ad}_{A_{if \nabla\bar{f}}} L_f = L_{  f  |\nabla f|^2}$.
\end{lem}
\noindent{\bf Remark:} In Lemma \ref{H0-4}, $g\cdot \nabla f$ denotes by extension the function which vanishes outside the support of $g$ and coincides with $g\cdot \nabla f$ inside (where $\nabla f$ is well-defined).

\subsection{Abstract results}

Let $f\in C^0({\mathbb T}^d)$, $|f|=1$ and consider a unitary operator $U$ defined on $l^2({\mathbb Z}^d)$. If $L_f^*U -I$ is compact (or equivalently if $U-L_f$ compact), we know from Weyl Theorem that $\sigma_{\text{ess}}(U)=\sigma_{\text{ess}}(L_f)=$ Ran $f$. Under some additional assumptions, we have that:
\begin{thm}\label{laurent-glob}{\bf (Global version)} Consider a non-constant symbol $f\in C^3({\mathbb T}^d)$ with $|f|=1$. Let $U\in {\cal C}^{1,1}(A_{if\nabla \bar{f}})$ be a unitary operator defined on $l^2({\mathbb Z}^d)$ such that $L_f^*U-I$ is compact. Then, $\sigma_{\text{ess}}(U)=$ Ran $f$ and
\begin{itemize}
\item[(a)] Given any Borel set $\Lambda\subset {\mathbb T}$ such that $e^{i\overline{\Lambda}}\subset$ Ran $f\setminus f(\kappa_f)$, $U$ has at most a finite number of eigenvalues in $e^{i\Lambda}$. Each of these eigenvalues has finite multiplicity.
\item[(b)] A LAP holds for $U$ on Ran $f\setminus \sigma_{\text{pp}}(U)\cup f(\kappa_f)$ w.r.t $A_{if\nabla \bar{f}}$ and $U$ has no singular continuous spectrum in Ran $f\setminus f(\kappa_f)$.
\end{itemize}
In particular, if $f(\kappa_f)$ has a finite number of accumulation points, $U$ has no singular continuous spectrum.
\end{thm}

\noindent{\bf Remark:} Consider a non-constant symbol $f\in C^3({\mathbb T}^d)$ with $|f|=1$ and let $U$ be defined as a multiplicative perturbation of $L_f$, $U=W_1L_fW_2$ where $W_1$ and $W_2$ are unitary operators defined on $l^2({\mathbb Z}^d)$. We already know that $L_f\in C^3(A_{if\nabla \bar{f}})\subset {\cal C}^{1,1}(A_{if\nabla \bar{f}})$ by Lemma \ref{H0-4}. So, if $W_1$ and $W_2$ belong to ${\cal C}^{1,1}(A_{if\nabla \bar{f}})$, then $U\in {\cal C}^{1,1}(A_{if\nabla \bar{f}})$ by Proposition \ref{4}.

The strategy behind the proof of Theorem \ref{laurent-glob} can also be applied to derive local spectral results. Before stating Theorem \ref{laurent-loc}, let us introduce the following language shortcut:
\begin{defin}\label{fgood} Let $f\in C^0({\mathbb T}^d)$ with $|f|=1$ and $\Lambda \subset {\mathbb T}$ be a Borel subset. We say that $\Lambda$ is $M_f$-good if
\begin{itemize}
\item $e^{i\overline{\Lambda}}\subset$ Ran $f\setminus f(\kappa_f)$,
\item there exists an open set $\Theta \subset {\mathbb T}^d$ on which $f$ is of class $C^3$ and $\overline{f^{-1}(e^{i\Lambda})}\subset \Theta$.
\end{itemize}
If $\eta$ is a smooth real-valued function compactly supported on $\Theta$ such that $\eta \upharpoonright \overline{f^{-1}(e^{i\Lambda})} \equiv 1$, we say that $\eta$ is $(M_f,\Lambda)$-adapted.
\end{defin}
\noindent{\bf Remark:} Let $f\in C^0({\mathbb T}^d)$ with $|f|=1$. If the Borel set $\Lambda \subset {\mathbb T}$ is $M_f$-good, then $\Lambda \subset$ Ran $f=\sigma_{ess}(L_f)$.

\begin{thm}\label{laurent-loc}{\bf (Local version)} Consider a non-constant symbol $f\in C^0({\mathbb T}^d)$ with $|f|=1$. Let $\Lambda$ be an open subset of ${\mathbb T}$, which is $M_f$-good and $\eta$ a $(M_f,\Lambda)$-adapted smooth real-valued function. Let $U\in {\cal C}^{1,1}(A_{i\eta f\nabla \bar{f}})$ be a unitary operator defined on $l^2({\mathbb Z}^d)$ such that $L_f^*U-I$ is compact. Then, $\sigma_{\text{ess}}(U)=$ Ran $f$ and
\begin{itemize}
\item[(a)] Given any Borel subset $\Lambda' \subset {\mathbb T}$ such that $\overline{\Lambda'}\subset \Lambda$, $U$ has at most a finite number of eigenvalues in $e^{i\Lambda'}$. Each of these eigenvalues has finite multiplicity.
\item[(b)] A LAP holds for $U$ on $e^{i\Lambda}\setminus \sigma_{\text{pp}}(U)$ w.r.t $A_{i\eta f\nabla \bar{f}}$ and $U$ has no singular continuous spectrum in $e^{i\Lambda}$.
\end{itemize}
\end{thm}
The proof of Theorems \ref{laurent-glob} and \ref{laurent-loc} are postponed to Sections 4.3 and 4.4 respectively. Let us point out that a generalized form of the LAP also holds for $U$ in Theorems \ref{laurent-glob} and \ref{laurent-loc} (see e.g. \cite{abc1}). The fine distribution properties of the point spectrum are not covered by our results.

The case of the translations on the lattice ${\mathbb Z}^d$ deserves a special treatment. If the function $f$ is defined on ${\mathbb T}^d$ by $f(\theta)= e^{i\alpha\cdot \theta}$ for some $\alpha\in {\mathbb Z}^d\setminus \{0\}$, then $L_f=T^{\alpha}$ and $A_{if\nabla \bar{f}} = \sum_{j=1}^d \alpha_j X_j$. Theorem \ref{laurent-glob} rewrites:
\begin{prop}\label{transl} Let $\alpha\in {\mathbb Z}^d\setminus \{0\}$. Let $U\in {\cal C}^{1,1}(\sum_{j=1}^d \alpha_j X_j)$ be a unitary operator defined on $l^2({\mathbb Z}^d)$ such that $\mathrm{ad}_{\sum_{j=1}^d \alpha_j X_j} T^{-\alpha} U$ is compact. Then,
\begin{itemize}
\item[(a)] $U$ has at most a finite number of eigenvalues in $\partial{\mathbb D}$ and each of these eigenvalues has finite multiplicity.
\item[(b)] A LAP holds for $U$ on $\partial{\mathbb D}\setminus \sigma_{pp}(U)$ w.r.t. $\sum_{j=1}^d \alpha_j X_j$ and $U$ has no singular continuous spectrum.
\end{itemize}
\end{prop}
The proof of Proposition \ref{transl} is postponed to Section 4.5.

\section{GGT Matrices}

GGT matrices appeared first in the theory of orthogonal polynomials on the unit circle \cite{g}. For an introduction to this subject in general and the model in particular, the reader is referred to \cite{S1}. The spectral analysis of such matrices has been undertaken in the contexts of periodic and random Verblunsky coefficients \cite{gt}, \cite{gnva1}, \cite{gnva2}, \cite{gol}, \cite{gn}, based on the theory of orthogonal polynomials and the associated transfer matrices formalism. In this section, we reconsider those GGT matrices with asymptotically constant Verblunsky coefficients and show how some existing spectral results can be recovered as a corollary of Theorem \ref{laurent-glob}.

Our description of the model follows \cite{gt}. In this section, $(e_k)_{k\in {\mathbb Z}}$ denotes the canonical orthonormal basis of $l^2({\mathbb Z})$. We denote by $T$,  $X$ the linear operators defined by:
\begin{eqnarray}
T e_k  & = & e_{k+1} \label{translation}  \\
X e_k & = & k e_k  \enspace . \label{Xope}
\end{eqnarray}
for all $k\in {\mathbb Z}$. If $\gamma=(\gamma_k)_{k\in {\mathbb Z}}$ is a sequence of complex numbers, we define the sequences $S\gamma$ and $\Delta\gamma$ by: $(S\gamma)_k = \gamma_{k+1}$ and $(\Delta\gamma)_k = \gamma_k-\gamma_{k+1}$ for all $k\in {\mathbb Z}$. We also write $\xi=(k)_{k\in {\mathbb Z}}$.

Consider a sequence $(\alpha_k)_{k\in {\mathbb Z}}\in {\mathbb D}^{\mathbb Z}$ such that:
\begin{equation*}
\sum_{k=0}^{\infty} |\alpha_k|^2 = \infty = \sum_{k=-1}^{-\infty} |\alpha_k|^2 \enspace ,
\end{equation*}
and define the sequence $(a_k)_{k\in {\mathbb Z}}\in [1,\infty)^{\mathbb Z}$ by: $a_k^{-2}+|\alpha_k|^2= 1$ for all $k\in {\mathbb Z}$. It follows from Lemma 2.2 \cite{gt} that the linear operator $H(\alpha)$ ('$H$' for Hessenberg) defined by:
\begin{eqnarray*}
H(\alpha) e_k = \frac{1}{a_k} \, e_{k-1} -\overline{\alpha_k} \sum_{l=k}^{\infty} \alpha_{l+1} \left(\prod_{m=k+1}^{l} \frac{1}{a_m}\right) e_l
\end{eqnarray*}
is unitary on $l^2({\mathbb Z})$. The operator $H(\alpha)$ is the GGT representation associated to the sequence of Verblunsky coefficients $(\alpha_k)$. For simplicity, we assume throughout this section that $0< \inf_{k\in {\mathbb Z}} |\alpha_k|\leq \sup_{k\in {\mathbb Z}} |\alpha_k| < 1$. Under this assumption, the operator can be rewritten as follows:
\begin{equation}\label{Uzero}
H(\alpha) = T^*D_2(\alpha) - T^*D_1(\alpha)T(I-D_2(\alpha)T)^{-1}D_1(\alpha)^* \, ,
\end{equation}
where $D_1(\alpha)$ and $D_2(\alpha)$ are the bounded diagonal operators defined on $l^2({\mathbb Z})$ by: $D_1(\alpha) e_k = \alpha_k e_k$ and $D_2(\alpha) e_k = a_k^{-1} e_k$ for all $k\in {\mathbb Z}$.

We start with the spectral analysis of such matrices, when the sequence of Verblunsky coefficients is a constant $\alpha_{\infty}$ with modulus different from $0$ and $1$ and define $a\in (1,\infty)$ via: $|\alpha_{\infty}|^2 + a^{-2}= 1$. The GGT representation associated to such a sequence is actually the Laurent operator associated to the function $f_a$:
\begin{equation}\label{Uzeroc}
L_{f_a} = \frac{1}{a}\, T^* - |\alpha |^2 \, \sum_{l=0}^{\infty}  \left( \frac{T}{a}\right)^l \enspace ,
\end{equation}
where the smooth complex-valued function $f_a$ is defined on ${\mathbb T}$ by:
\begin{equation}\label{falpha}
f_a(\theta) = - \frac{e^{-i\theta}-a}{e^{i\theta}-a}\enspace .
\end{equation}
$L_{f_a}$ is clearly purely absolutely continuous. We also have that:
\begin{equation*}
i f_a \bar{f_a}'(\theta) = - \frac{2-2a\cos \theta}{|1-ae^{i\theta}|^2}=2\Re (\frac{1}{ae^{i\theta}-1}) = \sum_{m\neq 0} a^{-|m|} e^{im\theta} \enspace .
\end{equation*}
for all $\theta \in {\mathbb T}$. In view of Section 2.3, we define the conjugate operator for $L_{f_a}$ by:
\begin{equation}\label{conjugate1}
A_a := A_{if_a \bar{f_a}'}= \frac{1}{2} \left(L_{i f_a \bar{f_a}'} X + XL_{i f_a \bar{f_a}'}\right) = \frac{1}{2} \sum_{m \neq 0} a^{-|m|} (T^m X+ XT^m)\, .
\end{equation}
According to Lemma \ref{sa}, ${\cal D}_X={\cal D}(X)$ is a core for $A_a$.

Now, we consider some local perturbations $H(\alpha)$ of $L_{f_a}$ through local fluctuations of the sequence $(\alpha_k)$. These fluctuations are measured by the family of norms $(q_n)_{n\geq 0}$ defined on ${\mathbb C}^{\mathbb Z}$ by:
\begin{eqnarray}\label{seminorm}
q_0(\gamma) & = & \| \gamma \|_{\infty}:=  \sup_{k \in {\mathbb Z}} |\gamma_k | \\
q_{n+1}(\gamma) & = & q_n(\gamma) + \|\xi^{n+1} \Delta^{n+1} \gamma \|_{\infty}  \enspace .
\end{eqnarray}

\begin{thm}\label{ggtperturbed1} Let $\alpha_{\infty}\in {\mathbb D}\setminus \{0\}$ and define $a\in (1,\infty)$ such that: $|\alpha_{\infty}|^2 +a^{-2}= 1$. Let $(\alpha_k) \in {\mathbb D}^{\mathbb Z}$ such that: $0<\inf_{k\in {\mathbb Z}} |\alpha_k|\leq \sup_{k\in {\mathbb Z}} |\alpha_k|< 1$. Assume that for all $k\in {\mathbb Z}$, $\alpha_k =\alpha_{\infty} (1+u_k+v_k+w_k)$ where $(u_k)$, $(v_k)$ and $(w_k)$ are complex valued sequences such that:
\begin{itemize}
\item 
\begin{equation*}
\int_1^{\infty} \sup_{b_1 r\leq |k| \leq b_2 r} |u_k|\, dr < \infty \enspace ,
\end{equation*}
for some $0<b_1<b_2< \infty$,
\item $q_{1}(v)<\infty$, $\lim_{|k|\rightarrow \infty} v_k=0$ and
\begin{equation*}
\int_1^{\infty} \sup_{b_1 r\leq |k|\leq b_2 r} |v_{k+1} -v_k|\, dr < \infty
\end{equation*}
for some $0<b_1<b_2<\infty$,
\item $q_{2}(w) < \infty$ and $\lim_{|k|\rightarrow \infty} w_k=0$.
\end{itemize}
Then, taking arguments in $[0,2\pi)$, we have that:
\begin{itemize}
\item $\sigma_{\text{ess}} (H(\alpha))=\sigma_{\text{ess}} (L_{f_a})=\Theta_a:=\{e^{i\theta}; \arg f_a(-\theta_a) \leq \theta \leq \arg f_a(\theta_a)\}$, with $\theta_a=\cos^{-1} (a^{-1})$,
\item any subset $\Theta$ of $\Theta_a$ whose closure does not contain any endpoints of $\Theta_a$ contains at most a finite number of eigenvalues. Each of these eigenvalues has finite multiplicity.
\item A LAP holds for $H(\alpha)$ on $\Theta_a \setminus \{f_a(\pm \theta_a)\}\cup \sigma_{\text{pp}}(H(\alpha))$ w.r.t $A_a$ and $H(\alpha)$ does not have any singular continuous spectrum.
\end{itemize}
\end{thm}

\noindent{\bf Remark:} The eigenvalues may accumulate at the endpoints of the arc $\Theta_a$. The reader is referred to \cite{gnva1}, \cite{gnva2} for further results.

The proof of Theorem \ref{ggtperturbed1} is carried out in three steps. We translate the properties of the sequence $\alpha$ in terms of the diagonal operator $D_1(\alpha)$ (Section 3.1) and then in terms of the unitary operator $H(\alpha)$ (Section 3.2). These results are combined with Theorem \ref{laurent-glob} in Section 3.3.

\subsection{From $\alpha$ to $D_1(\alpha)$}

Let us adopt some local notations. To any bounded sequence $\gamma:=(\gamma_k)_{k\in {\mathbb Z}}$ in ${\mathbb C}^{\mathbb Z}$, we associate the bounded linear operator $D_{\gamma}$ defined by its action on the canonical orthonormal basis of $l^2({\mathbb Z})$: $D_{\gamma} e_k=\gamma_k e_k$, $k\in {\mathbb Z}$. The operator $D_{\gamma}$ is normal and we recall that $\|D_{\gamma}\|=\sup_{k}|\gamma_k|=q_0(\gamma)$. Note that if $\gamma$ and $\beta$ are two bounded sequences in ${\mathbb C}^{\mathbb Z}$ and $c\in {\mathbb C}$, then $D_{\gamma+\beta}=D_{\gamma}+D_{\beta}$, $D_{\gamma\cdot \beta}=D_{\gamma}D_{\beta}$, $D_{c\gamma}=cD_{\gamma}$ and $[D_{\gamma},D_{\beta}]=0$. In addition, if the sequence $\xi\cdot\gamma$ is bounded, then $D_{\xi\cdot \gamma}=XD_{\gamma}=D_{\gamma}X$.

For any bounded sequence $\gamma$, $D_{\gamma} \in C^{\infty}(X)$ and $\mathrm{ad}_X D_{\gamma}=0$. Along the next lines, we also apply the following identities without further comments: for all $m\in {\mathbb Z}$ and any bounded sequence $\gamma$ in ${\mathbb C}^{\mathbb Z}$, $T^m D_{\gamma}= D_{S^{-m}\gamma}T^m$, $D_{\gamma}T^m=T^mD_{S^m \gamma}$.
\begin{lem}\label{q0q1q2} Let $\gamma$ be a bounded sequence of ${\mathbb C}^{\mathbb Z}$.
\begin{itemize}
\item[(a)] If $q_1(\gamma)< \infty$, then there exists $C>0$ such that $q_0(\xi(\gamma -S^m \gamma))\leq C m^2 q_1(\gamma)$, for all $m\in {\mathbb Z}$,
\item[(b)] If $q_2(\gamma)< \infty$, then there exists $C>0$ such that $q_1(\xi(\gamma -S^m \gamma))\leq C |m|^3 q_2(\gamma)$, for all $m\in {\mathbb Z}$.
\end{itemize}
\end{lem}
\noindent{\bf Proof:} We recall that if $\gamma$ is a bounded sequence, $\|S^m\gamma\|_{\infty}=\|\gamma\|_{\infty}$ for any $m\in {\mathbb Z}$ and $\|\Delta\gamma \|_{\infty}\leq 2\|\gamma\|_{\infty}$. If $\xi^p \gamma$ is a bounded sequence for some $p\in {\mathbb N}$, then $\|\xi^l \gamma \|_{\infty}\leq \|\xi^p \gamma\|_{\infty}$ for any $l\in \{1,\ldots,p\}$. We restrict our discussion to the case $m>0$ (the case $m<0$ can be carried out similarly):
$$
\xi (\gamma - S^m\gamma) =\sum_{l=0}^{m-1} S^l(\xi-l) \Delta\gamma
$$
which proves statement (a), and
\begin{equation*}
\xi\Delta(\xi (\gamma - S^m\gamma)) = \sum_{l=0}^{m-1} S^l(\xi-l)^2 \Delta^2\gamma-\sum_{l=1}^m S^l(\xi-l)\Delta\gamma  \enspace ,
\end{equation*}
which allows to get (b). \ep

\begin{prop}\label{criteria1} Let $a\in (1,\infty)$. If $q_1(\gamma)<\infty$, then $D_{\gamma}\in C^1(A_a)$.
\end{prop}
\noindent{\bf Proof:} We recall that: $[T^m, D_{\gamma}] = T^m D_{(\gamma - S^m\gamma)}=D_{(S^{-m}\gamma-\gamma)} T^m$ for all $m\in {\mathbb N}$. Assume first that $q_1(\gamma)<\infty$. As a sesquilinear form on ${\cal D}(X)\times {\cal D}(X)$, one has that
\begin{equation*}
A_a D_{\gamma} - D_{\gamma} A_a = \frac{1}{2} \sum_{m\neq 0} a^{-|m|} \left([T^m,D_{\gamma}]X+X[T^m,D_{\gamma}] \right) = \frac{1}{2}\sum_{m\neq 0} a^{-|m|}  \left(T^m D_{(\gamma-S^m \gamma)\xi}-D_{(\gamma-S^{-m} \gamma)\xi} T^m  \right) \enspace .
\end{equation*}
Using Lemma \ref{q0q1q2} and the fact that $T$ is unitary, the RHS of the previous identity defines a norm convergent series of bounded operators ($|a|>1$). Since ${\cal D}(X)$ is a core for $A_a$, this implies that $D_{\gamma}\in C^1(A_a)$ with
\begin{equation*}
\mathrm{ad}_{A_a} D_{\gamma} = \frac{1}{2}\sum_{m\neq 0} a^{-|m|}  \left(T^m D_{(\gamma-S^m \gamma)\xi}-D_{(\gamma-S^{-m} \gamma)\xi} T^m  \right) \enspace ,
\end{equation*}
and that $\|\mathrm{ad}_{A_a} D_{\gamma}\| \leq C_a q_1(\gamma)$ for some positive constant $C_a$ which depends only on $a$. \ep

\begin{prop}\label{criteria2} Let $a\in (1,\infty)$. If $q_2(\gamma)<\infty$, then $D_{\gamma}\in C^2(A_a)$.
\end{prop}
\noindent{\bf Proof:} If $q_2(\gamma)<\infty$, then $q_1(\gamma)<\infty$ and $D_{\gamma}\in C^1(A_a)$ by Proposition \ref{criteria1}. It remains to prove that $\mathrm{ad}_{A_a} D_{\gamma}$ belongs to $C^1(A_a)$. By Lemma \ref{H0-4}, we have that for all $m\in {\mathbb Z}$, $T^m \in C^{\infty}(A_a)$, $\mathrm{ad}_{A_a} T^m = m T^m L_{if_a \bar{f_a'}}$ (recall that $T$ commutes with $L_{if_a \bar{f_a'}}$) hence $\| \mathrm{ad}_{A_a} T^m \| \leq C_a |m|$ for some $C_a>0$. By Lemma \ref{q0q1q2} and Proposition \ref{criteria1}, we have also that for all $m\in {\mathbb Z}$, $q_1((\gamma-S^m \gamma)\xi) < \infty$, $D_{(\gamma-S^m \gamma)\xi}\in C^1(A_a)$ and $\|\mathrm{ad}_{A_a} D_{(\gamma-S^m \gamma)\xi}\| \leq C_a |m|^3 q_2(\gamma)$ for some $C_a>0$. It follows that for all $m\in {\mathbb Z}$, the operators $T^m D_{(\gamma-S^m \gamma)\xi}-D_{(\gamma-S^{-m} \gamma)\xi} T^m \in C^1(A_a)$ (see e.g. Proposition \ref{3}) and that the series
\begin{equation*}
\sum_{m \neq 0} a^{-|m|}  \mathrm{ad}_{A_a}(T^m D_{(\gamma-S^m \gamma)\xi}-D_{(\gamma-S^{-m} \gamma)\xi} T^m)
\end{equation*}
is norm convergent. The conclusion follows from Lemma \ref{H0-1}. \ep

The next result provides a practical criterion to deal with fractional order regularities. See e.g. Theorem 6.1 in \cite{sah} for a proof.
\begin{thm}\label{gsah} Let $s\in \{0,1\}$ and $Q$ be a self-adjoint operator in ${\cal H}$ bounded from below by a strictly positive constant such that $A^l Q^{-l}$ is continuous for some integer $l\in {\mathbb N}$, $s < l$. Then, a bounded symmetric operator $B$ is of class ${\cal C}^{s,1}(A)$ if there exists a function $\chi \in C^{\infty}_0({\mathbb R})$ which is positive on some interval $(b_1,b_2)$, $0<b_1<b_2<\infty$, such that:
\begin{equation}\label{gsah2}
\int_1^{\infty} \|r^s\chi(Q/r)B\| \frac{dr}{r} < \infty
\end{equation}
\end{thm}
In the following, the roles of $A$ and $Q$ are endorsed by $A_a$ and $\bra X\ket$ respectively (see Lemma \ref{sa}). In order to obtain estimates (\ref{gsah2}), it is enough to prove that:
\begin{equation*}
\int_R^{\infty} \|r^s\chi(Q/r)B\| \frac{dr}{r} < \infty
\end{equation*}
for some $R>0$.

\begin{prop}\label{d1approx} Let $\gamma$ be a bounded sequence of complex numbers.
\begin{itemize}
\item[(a)] If
\begin{equation}\label{h1}
\int_1^{\infty} \sup_{b_1 r\leq |k|\leq b_2 r} |\gamma_k|\, dr < \infty
\end{equation}
for some $0<b_1 < b_2 < \infty$, then $D_{\gamma} \in {\cal C}^{1,1}(A_a)$.
\item[(b)] If $q_1(\gamma)<\infty$ and
\begin{equation}\label{h2}
\int_1^{\infty} \sup_{b_1 r\leq |k|\leq b_2 r} |\gamma_{k+1} -\gamma_k|\, dr < \infty
\end{equation}
for some $0<b_1<b_2<\infty$, then $D_{\gamma} \in C^1(A_a)$ and $\mathrm{ad}_{A_a} D_{\gamma} \in {\cal C}^{0,1}(A_a)$. In particular, $D_{\gamma} \in {\cal C}^{1,1}(A_a)$.
\end{itemize}
\end{prop}
\noindent{\bf Proof:} We first establish Proposition \ref{d1approx} for bounded sequences of real numbers and then extend it to the complex case:\\
{\bf Real case.} We start with statement (a) by assuming (\ref{h1}). Let $\chi$ be any smoothed characteristic function supported on $(b_1+\epsilon,b_2-\epsilon)$ with $0<\epsilon<(b_2-b_1)/2$. We have that for $R>1$ big enough,
\begin{equation*}
\int_R^{\infty} \|\chi(\bra X\ket/r)D_{\gamma}\| \, dr \leq  \int_R^{\infty} \sup_{b_1 r\leq |k|\leq b_2 r} |\gamma_k|\, dr < \infty \enspace ,
\end{equation*}
and the conclusion follows from Theorem \ref{gsah}. Now, we prove statement (b) assuming (\ref{h2}). According to Proposition \ref{criteria1} , if $q_1(\gamma) < \infty$, then $D_{\gamma}\in C^1(A_a)$ and we can rewrite:
\begin{equation*}
\mathrm{ad}_{A_a} D_{\gamma} = \frac{1}{2}\left( \sum_{m\neq 0} a^{-|m|}  T^m D_{(\gamma-S^m \gamma)\xi}- \sum_{m\neq 0} a^{-|m|} D_{(\gamma-S^m \gamma)\xi} T^{-m} \right) \enspace .
\end{equation*}
It remains to prove that $\mathrm{ad}_{A_a} D_{\gamma} \in {\cal C}^{0,1}(A_a)$. ${\cal C}^{0,1}(A_a)$ is a linear subspace which is stable under adjunction $^*$ (see e.g. Proposition \ref{4}). Since $\gamma$ is real-valued, we observe that $(\sum_{m\neq 0} a^{-|m|}T^m D_{(\gamma-S^m \gamma)\xi})^*=\sum_{m\neq 0} a^{-|m|}D_{(\gamma-S^m \gamma)\xi} T^{-m}$. So, it is enough to prove that the norm convergent series
$$S:=\sum_{m\neq 0} a^{-|m|} D_{(\gamma-S^m \gamma)\xi} T^{-m}$$
belongs to ${\cal C}^{0,1}(A_a)$. Let $\chi$ be any smoothed characteristic function supported on $(b_1+2\epsilon,b_2-2\epsilon)$ with $0<\epsilon<(b_2-b_1)/4$. In view of Theorem \ref{gsah}, we will prove that the map $r\mapsto r^{-1} \|\chi(\bra X\ket/r) S\|$ is integrable on $[R,\infty)$ for some $R>1$. For any $r\geq 1$, we define $N_r=[\epsilon\ln r]$. We observe that for all $r>1$, $\|\chi(\bra X\ket/r) S\| \leq S_1(r)+S_2(r)$ where
\begin{eqnarray*}
S_1(r) &=& \sum_{|m|< N_r } a^{-|m|} \| \chi(\bra X\ket /r) D_{\xi(\gamma -S^m \gamma)} \|\\
S_2(r) &=& \sum_{|m|\geq N_r } a^{-|m|} \| \chi(\bra X\ket /r) D_{\xi(\gamma -S^m \gamma)} \|
\end{eqnarray*}
By Lemma \ref{q0q1q2}, there exists $C>0$ such that for all $r>1,$
\begin{equation*}
S_2(r) \leq C \sum_{|m|\geq N_r } a^{-|m|} \sup_{k\in {\mathbb Z}} |k(\gamma_k-\gamma_{k+m})| \leq C \sum_{|m|\geq N_r } a^{-|m|} m^2 q_1(\gamma) \leq C\frac{\ln^2 r}{r^{\beta}}\enspace ,
\end{equation*}
where $\beta=\ln |a|$, $|a|>1$, which implies the integrability of the map $r\mapsto r^{-1}S_2(r)$. Now, for $r\geq R$, $R>1$ large enough, we have that:
$$
S_1(r) \leq \sum_{|m|< N_r} a^{-|m|} \sup_{(b_1+\epsilon)r \leq |k|\leq (b_2-\epsilon)r} |k(\gamma_{k+m}-\gamma_k)| \leq b_2 r \sum_{|m|< N_r} a^{-|m|} \sup_{(b_1+\epsilon)r \leq |k|\leq (b_2-\epsilon)r} |\gamma_{k+m}-\gamma_k| \enspace .
$$
Combining the facts that $|m|\leq \epsilon \ln r< \epsilon r$ and $(b_1+\epsilon)r \leq |k|\leq (b_2-\epsilon)r$, we have that: $\sup_{(b_1+\epsilon)r \leq |k|\leq (b_2-\epsilon)r}  |\gamma_k - \gamma_{k+m} | \leq |m| \sup_{b_1 r \leq |k|\leq b_2 r} |\gamma_k - \gamma_{k+1}|$ for all $|m|\leq N_r$. We deduce that for all $r\geq R,$
$$
S_1(r)\leq  b_2 r \left(\sum_{m\in {\mathbb Z}} |m| a^{-|m|}\right) \sup_{b_1 r\leq |k|\leq b_2 r} |\gamma_k - \gamma_{k+1}| \enspace .
$$
Due to the hypotheses, this implies the integrability of the map $r\mapsto r^{-1} S_1(r)$ and concludes the discussion for bounded real sequences. \\
{\bf Complex case.} Now, assume that $\gamma$ is a bounded sequence of complex numbers which satisfies (\ref{h1}). The integrability of the map $r\mapsto \sup_{b_1 r \leq |k|\leq b_2 r} |\gamma_k|$ is equivalent to the joint integrability of $r\mapsto \sup_{b_1 r\leq |k|\leq b_2 r} |Re \gamma_k|$ and $r\mapsto \sup_{b_1 r\leq |k|\leq b_2 r} |\Im \gamma_k|$. It follows from the real case that $D_{\Re \gamma}$ and $D_{\Im \gamma}$ belongs to ${\cal C}^{1,1}(A_a)$. Since ${\cal C}^{1,1}(A_a)$ is a linear space and $D_{\gamma}=D_{\Re \gamma}+iD_{\Im \gamma}$, statement (a) follows. The proof of statement (b) is similar. \ep

In the next section, we show how the properties of the diagonal and normal operator $D_1(\gamma)$ are related to those of the unitary operator $H(\gamma)$.

\subsection{From $D_1(\alpha)$ to $H(\alpha)$}

In this section, $\gamma$ and $\beta$ stand for two generic sequences of ${\mathbb D}^{\mathbb Z}$ such that: $0<b:=\inf_k (|\beta_k|^2,|\gamma_k|^2)$ and $B=\sup_k (|\beta_k|^2,|\gamma_k|^2)<1$. Let $\epsilon\in (0,1-B)$. There exists $\Phi\in C_0^{\infty}({\mathbb R}, {\mathbb R})$ with compact support in $(-\infty,1-\epsilon)$ such that:
\begin{equation*}
D_2(\gamma) = \sqrt{1-|D_1(\gamma)|^2} =\Phi(|D_1(\gamma)|^2) \quad \text{and} \quad D_2(\beta) =\Phi(|D_1(\beta)|^2)
\end{equation*}
where $|D_1(\gamma)|^2=D_1(\gamma)^*D_1(\gamma)$ and $|D_1(\beta)|^2=D_1(\beta)^*D_1(\beta)$ (the operators $D_1(\gamma)$ and $D_1(\beta)$ are normal). It follows that:
\begin{lem}\label{ggt0} If $D_1(\gamma)-D_1(\beta)$ is compact, so are $D_2(\gamma)-D_2(\beta)$ and $H(\gamma)-H(\beta)$.
\end{lem}
\noindent{\bf Proof:} The operators $D_1(\gamma)^*-D_1(\beta)^*$ and $|D_1(\gamma)|^2-|D_1(\beta)|^2$ are compact. The first statement follows from the Stone-Weierstrass Theorem. We have that:
\begin{align}\label{ha-hb}
H(\gamma)-H(\beta) &= T(D_2(\gamma)-D_2(\beta)) - T^*(D_1(\gamma)-D_1(\beta))T (1-D_2(\gamma)T)^{-1} D_1(\gamma)^* \nonumber \\
& + T^*D_1(\beta))T (1-D_2(\gamma)T)^{-1} (D_2(\gamma)-D_2(\beta))T  (1-D_2(\gamma)T)^{-1}D_1(\gamma)^* \nonumber \\
& - T^*D_1(\beta) T (1-D_2(\beta)T)^{-1} (D_1(\gamma) - D_1(\beta))^* \enspace .
\end{align}
Each term on the RHS is the product of bounded operators with at least one compact operator, which concludes the proof. \ep

In other words, the difference $H(\gamma)-H(\beta)$ is compact whenever:
\begin{equation*}
\lim_{|k|\rightarrow \infty} (\gamma_k -\beta_k) = 0
\end{equation*}
In this case, we conclude from Weyl's Theorem that: $\sigma_{\text{ess}}(H(\gamma)) = \sigma_{\text{ess}}(H(\beta))$.

\begin{lem}\label{ggt2} Let $n\in {\mathbb N}$ and $(s,p)\in [0,\infty)\times [1,\infty)$. If $D_1(\gamma)$ belongs to $C^n(A_a)$ (resp. ${\cal C}^{s,p}(A_a)$), then $D_2(\gamma)$ and $H(\gamma)$ also belong to $C^n(A_a)$ (resp. ${\cal C}^{s,p}(A_a)$).
\end{lem}
\noindent{\bf Proof:} The operators $D_1(\gamma)^*$ and $|D_1(\gamma)|^2$ also belong to $C^n(A_a)$ (resp. ${\cal C}^{s,p}(A_a)$). The fact that $D_2(\gamma)$ belongs to $C^n(A_a)$ (resp. ${\cal C}^{s,p}(A_a)$) follows from Theorem 6.2.5 and Corollary 6.2.6 in \cite{abmg}. According to Lemma \ref{H0-4}, $T$ and $T^*$ belong to $C^{\infty}(A_a)$. Therefore, the last part of the proof follows from formula (\ref{Uzero}) by using the properties of the classes $C^n(A_a)$ and ${\cal C}^{s,p}(A_a)$ (see e.g. Propositions \ref{3} and \ref{4}). \ep

The proof of Lemma \ref{ggt2} can also be carried out by means of the Helffer-Sj\"ostrand functional calculus.
 
\subsection{Proof of Theorem \ref{ggtperturbed1}}

The function $f_a$ belongs to $C^{\infty}({\mathbb T})$ and $|f_a|=1$. We also observe that $D_1(\alpha)-\alpha_{\infty}I$ is compact. Gathering Propositions \ref{criteria2}, \ref{d1approx}, Lemmata \ref{ggt0}, \ref{ggt2}, we have that $H(\alpha)\in {\cal C}^{1,1}(A_a)$ and that $L_{f_a}^* H(\alpha)-I=L_{f_a}^*(H(\alpha)-L_{f_a})$ is a compact operator. The conclusions follow by applying Theorem \ref{laurent-glob}, once noted that $f_a({\mathbb T})=$ Ran $f_a=\Theta_a$ and $f_a(\kappa_{f_a})=\{f_a(\pm\theta_a)\}$.

\section{Technicalities}

The proofs of Theorems \ref{laurent-glob}, \ref{laurent-loc} and Proposition \ref{transl} rely on a suitable application of the regular Mourre Theory, which is summed up in Section 4.1. We recall that to any bounded Borel function $\Phi$ on $\partial {\mathbb D}$ is associated a unique function $\phi$ defined on ${\mathbb T}$ by:  $\phi(\theta)=\Phi(e^{i\theta})$, for all $\theta \in {\mathbb T}$.

\subsection{Regular Mourre Theory for Unitary Operators}

We refer the reader to \cite{abc1} for the proofs of the abstract results presented here. Throughout this section, ${\cal H}$ is an infinite dimensional Hilbert space and $A$ denotes a self-adjoint operator (densely) defined on ${\cal H}$ with domain ${\cal D}(A)$. If $U$ is a unitary operator acting on ${\cal H}$, its spectral measure is denoted by $(E_{\Delta}(U))_{\Delta \in {\cal B}({\mathbb T})}$.

We use the following notations: if $(S,T)\in {\cal B}({\cal H})\times {\cal B}({\cal H})$, we write $S\simeq T$ if $S-T$ is a compact operator and $S\lesssim T$ (resp. $S\gtrsim T$) if $S\leq T+K$ (resp. $S\geq T+K$) for some compact operator $K$.

Let us start with a reformulation of Lemma 2.1 in \cite{abc1}. We note that a unitary operator $U$ belongs to $C^1(A)$ iff $U^*=U^{-1}$ belongs to $C^1(A)$. Then, for all $n\in {\mathbb Z}$, $U^n\in C^1(A)$ and $U^n {\cal D}(A)={\cal D}(A)$ (see e.g. Propositions \ref{2} and \ref{3}). 

\begin{lem}\label{propag0} Let $U$ be a unitary operator defined on ${\cal H}$. $U$ (or equivalently $U^*$) belongs to $C^1(A)$ if and only if one of the following statements is satisfied.
\begin{itemize}
\item[(a)] $U {\cal D}(A)\subset {\cal D}(A)$ and the sesquilinear form $F_{+}$ defined on ${\cal D}(A)\times {\cal D}(A)$ by $F_{+}(\varphi,\psi):=\langle U\varphi, AU\psi\rangle\ - \langle\varphi,A\psi\rangle$ is continuous for the topology induced by ${\cal H}\times {\cal H}$.
\item[(b)] There exists a core for $A$, denoted ${\cal S}$, such that $U {\cal S}\subset {\cal S}$ and the sesquilinear form $G_+$ defined on ${\cal S}\times {\cal S}$ by $G_+(\varphi,\psi):=\langle U\varphi, AU\psi\rangle\ - \langle\varphi,A\psi\rangle$ is continuous for the topology induced by ${\cal H}\times {\cal H}$.
\item[(c)] $U^* {\cal D}(A)\subset {\cal D}(A)$ and the sesquilinear form $F_{-}$ defined on ${\cal D}(A)\times {\cal D}(A)$ by $F_{-}(\varphi,\psi):=\langle \varphi,A\psi\rangle\ - \langle U^*\varphi,AU^*\psi\rangle$ is continuous for the topology induced by ${\cal H}\times {\cal H}$.
\item[(d)] There exists a core for $A$, denoted ${\cal S}$, such that $U^* {\cal S}\subset {\cal S}$ and the sesquilinear form $G_-$ defined on ${\cal S}\times {\cal S}$ by $G_-(\varphi,\psi):=\langle \varphi, A\psi\rangle\ - \langle U^*\varphi,AU^*\psi\rangle$ is continuous for the topology induced by ${\cal H}\times {\cal H}$.
\end{itemize}
If $(U^*AU-A)^o$ (resp. $B_+$, resp. $(A-UAU^*)^o$, resp. $B_-$) denotes the (unique) bounded operator associated to the continuous extension of $F_{+}$ (resp. $G_+$, resp. $F_{-}$, resp. $G_-$) to ${\cal H}\times {\cal H}$, then we have that: $(U^*AU-A)^o=B_+=U^*\mathrm{ad}_A U$ and $(A-UAU^*)^o=B_-=-U\mathrm{ad}_A U^*$.
\end{lem}
\noindent{\bf Proof:} Statement (a) implies obviously (b), while (c) implies (d) (by taking ${\cal S}={\cal D}(A)$). $U\in C^1(A)$ implies (a) and (c) by Lemma \ref{propag01}. By Lemma \ref{propag02}, (b) implies that $U\in C^1(A)$. By Lemma \ref{propag02bis}, (d) implies $U^*\in C^1(A)$, hence the result. \ep

For more details, we refer to Section 6.2. In the following, we drop the superscript $^o$ and write $(U^* AU-A)$ (resp. $(A-UAU^*)$) to refer to the bounded operator associated to the continuous extension of the sesquilinear forms $F_{+}$, $G_+$ (resp. $F_{-}$, $G_-$).

\noindent{\bf Remark:} If $U\in C^1(A)$ is a unitary operator, $U^*U=I=UU^*$ hence $U^*(\mathrm{ad}_A U)+(\mathrm{ad}_A U^*)U=0=U(\mathrm{ad}_A U^*)+(\mathrm{ad}_A U)U^*$. We deduce from Lemma \ref{propag0} that $(U^*AU-A)= U^*(A-UAU^*) U$, $(U^*AU-A)= U^*(\mathrm{ad}_A U)= -(\mathrm{ad}_A U^*)U$ and that $(A-UAU^*)= -U(\mathrm{ad}_A U^*)=(\mathrm{ad}_A U)U^*$. We will use these identities without any further comments.

Next we introduce the concept of Mourre estimates in the context of unitary operators:
\begin{defin}\label{propagating} Let $U$ be a unitary operator which belongs to $C^1(A)$. For a given $\Theta \in {\cal B}({\mathbb T})$, we say that
\begin{itemize}
\item $U$ is propagating w.r.t. $A$ on ${\Theta}$ if there exist $c >0$ such that: $E_{\Theta}(U) (U^*AU-A) E_{\Theta}(U) \gtrsim c E_{\Theta}(U)$
\item $U$ is strictly propagating w.r.t. $A$ on ${\Theta}$ if there exist $c >0$ such that: $E_{\Theta}(U) (U^*AU-A) E_{\Theta}(U) \geq c E_{\Theta}(U)$.
\item $U$ is weakly propagating w.r.t. $A$ if $(U^*AU-A)>0$, i.e. $\bra\varphi,(U^*AU-A)\varphi\ket >0$ for all $\varphi\in {\cal H}\setminus \{0\}$.
\end{itemize}
\end{defin}

\noindent{\bf Remark:} Since the spectral projectors associated to $U$ commute with $U$ and $U^*$, Definition \ref{propagating} can be equivalently stated by substituting the operator $(A-UAU^*)$ to $(U^*AU-A)$.

Now, let us state the main conclusions of the regular Mourre Theory for unitary operators. We start with the Virial Theorem (see e.g. Theorem 5.1 in \cite{abc1} for a proof):
\begin{thm} Let $U$ be in $C^1(A)$. Then, $E_{\{\theta \}}(U)(U^* AU-A)E_{\{\theta \}}(U)=0$ for all $\theta\in {\mathbb T}$.
\end{thm}
As immediate consequences, we have that:
\begin{itemize}
\item If $U$ is strictly propagating w.r.t $A$ on some Borel set $\Theta\subset {\mathbb T}$, then $U$ has no eigenvalue in $e^{i\Theta}$.
\item If $U$ is weakly propagating w.r.t $A$, then $U$ has no eigenvalue.
\end{itemize}
The next two results are also proved in \cite{abc1} (Corollary 4.1 and Theorem 2.3 respectively):
\begin{cor}\label{virial-2} Assume that $U$ is propagating w.r.t $A$ on the Borel subset $\Theta\subset {\mathbb T}$. Then, $U$ has at most a finite number of eigenvalues in $e^{i\Theta}$. Each of these eigenvalues has finite multiplicity.
\end{cor}
\begin{thm}\label{nosc2} Let $\Theta\subset {\mathbb T}$ be an open set. Assume that $U$ is propagating w.r.t $A$ on $\Theta$ and $U\in {\cal C}^{1,1}(A)$. Then, a LAP holds for $U$ on $\Theta \setminus \sigma_{\text{pp}}(U)$ w.r.t $A$ and $U$ has no singular continuous spectrum in $\Theta$.
\end{thm}

As a final remark, we observe that Mourre estimates are preserved under compact perturbations:
\begin{lem}\label{prop-compact} Let $U$ and $V$ be two unitary operators defined on ${\cal H}$, which belong to $C^1(A)$, $(c,C)\in {\mathbb R}^2$ and $B\in {\cal B}({\cal H})$.
\begin{itemize}
\item If $\mathrm{ad}_A U^*V\simeq 0$ and $B$ commutes with $U^*V$, then $(U^*AU-A)\simeq B$ if and only if $(V^*AV-A) \simeq B$.
\item If $U^*V\simeq I$ and $\mathrm{ad}_A (U^*V)\simeq 0$, then $(U^*AU-A)\simeq (V^*AV-A)$. In addition, for any real-valued $\Phi\in C^0(\sigma(U)\cup\sigma(V))$ (e.g. $\phi \in C^0({\mathbb T})$),
\begin{eqnarray*}
\Phi(U) (U^*AU-A) \Phi(U)\gtrsim c \Phi(U)^2 & \text{iff} & \Phi(V) (V^*AV-A) \Phi(V)\gtrsim c \Phi(V)^2\\
\Phi(U) (U^*AU-A) \Phi(U)\lesssim C \Phi(U)^2 & \text{iff} & \Phi(V) (V^*AV-A) \Phi(V)\lesssim C \Phi(V)^2
\end{eqnarray*}
\end{itemize}
\end{lem}
\noindent{\bf Proof:} Since $U$ and $V$ belong to $C^1(A)$, $W:= U^*V$ belongs to $C^1(A)$ (see e.g. Proposition \ref{3}). We can split the commutator $(V^* AV-A)$ as follows ($V=UW$):
\begin{eqnarray}
(V^* AV-A) &=& W^* (U^*AU-A)W+ (W^*AW-A) \label{vava1}\\
&=& (U^*AU-A) + \left(W^* (U^*AU-A)W-(U^*AU-A)\right) + (W^*AW-A) \label{vava2} \enspace ,
\end{eqnarray}
where the terms between parentheses are bounded. Now, we prove the first statement. Note that $\mathrm{ad}_A U^*V \simeq 0$ iff $\mathrm{ad}_A V^*U \simeq 0$ since $\mathrm{ad}_A V^*U=\mathrm{ad}_A (U^*V)^*=- (\mathrm{ad}_A U^*V)^*$. Since $U^*V$ is unitary, $U^*V$ commutes with $B$ iff $V^*U$ commutes with $B$. Indeed, since $W^*W=I=WW^*$, we have that $[B,W]W^*+W[B,W^*]=0$, which implies our claim. So, due to the symmetry of the problem, it is enough to prove one implication. The conclusion follows directly from identity (\ref{vava1}) and the fact that $(W^*AW-A)=W^* \mathrm{ad}_{A} W$. Let us prove the second statement. By hypothesis, we have that the terms $W^* (U^*AU-A)W-(U^*AU-A)=W^* \mathrm{ad}_{(U^*AU-A)} W=W^* \mathrm{ad}_{(U^*AU-A)} (W-I)$ and $(W^*AW-A)=W^* \mathrm{ad}_{A} W$ are compact. So $(U^*AU-A)\simeq (V^*AV-A)$ in view of identity (\ref{vava2}). Turning to the proof of the last equivalences, we note that due to the symmetry of the problem ($U^*V\simeq I$ iff $V^*U\simeq I$), it is again enough to prove one implication. We already know that for any real-valued $\Phi\in C^0(\sigma(U)\cup\sigma(V))$, the operator $\Phi(V) (V^* AV-A) \Phi(V) \simeq \Phi(V) (U^*AU-A) \Phi(V)$. Since $U \simeq V$, $\Phi(U)\simeq \Phi(V) $ (by Stone-Weierstrass Theorem). Combined with the fact that $U$ belongs to $C^1(A)$ this implies that $\Phi(V) (U^*AU-A) \Phi(V) \simeq \Phi(U) (U^*AU-A) \Phi(U)$. So far, we have proven that $\Phi(V) (V^*AV-A) \Phi(V) \simeq \Phi(U) (U^*AU-A) \Phi(U)$. The conclusion follows once noted that $\Phi(V)^2 \simeq \Phi(U)^2$ since $\Phi^2\in C^0(\sigma(U)\cup\sigma(V))$. \ep

\noindent{\bf Remark:} Since $U^*V-I=U^*(V-U)$, we have that $U^*V\simeq I$ iff $U\simeq V$. If the operators $U$ and $V$ belong to $C^1(A)$ and if $U^*V \simeq I$, then $\mathrm{ad}_A U^*V \simeq 0$ iff $\mathrm{ad}_A (V-U)\simeq 0$. This follows from the fact that: $\mathrm{ad}_A (U^*V) = \mathrm{ad}_A (U^*V-I) = \left(\mathrm{ad}_A U^*\right) (V-U)+ U^* \mathrm{ad}_A (V-U)$ and $\mathrm{ad}_A (U-V) = \left(\mathrm{ad}_A U\right) (I-U^*V)- U \mathrm{ad}_A (U^*V)$.

\subsection{Mourre estimates}

In this section, we establish some Mourre estimates for unitary Laurent operators (Proposition \ref{positivity}) and for suitable perturbations of them (Propositions \ref{positivity2} and \ref{positivity3}).

\noindent{\bf Notations.} Let $f\in C^0({\mathbb T}^d)$ with values in $\partial {\mathbb D}$ and $g\in C^0({\mathbb T}^d)$ be a real-valued function. Given any Borel set $\Lambda \subset {\mathbb T}$ such that $e^{i\Lambda} \cap$ Ran $f \neq \emptyset$, we define:
\begin{equation*}
c_{\Lambda,f,g}:=\min_{\theta \in \overline{f^{-1}(e^{i\Lambda})}}g(\theta) \quad, \quad C_{\Lambda,f,g}:=\max_{\theta \in \overline{f^{-1}(e^{i\Lambda})}}g(\theta) \enspace .
\end{equation*}
If ${\cal I}_{\Lambda}$ denotes the collection of all open sets $\Lambda' \subset {\mathbb T}$ such that $\overline{\Lambda}\subset \Lambda'$, we also write:
\begin{equation*}
c_{\Lambda,f,g}^{\sharp}:=\sup_{\Lambda'\in {\cal I}_{\Lambda}} c_{\Lambda',f,g} \quad, \quad C_{\Lambda,f,g}^{\flat}:=\inf_{\Lambda'\in {\cal I}_{\Lambda}} C_{\Lambda',f,g} \enspace .
\end{equation*}
\noindent{\bf Remark:} We note that $\overline{f^{-1}(e^{i\Lambda})}$ is a compact subset of ${\mathbb T}^d$. Clearly, $c_{\Lambda,f,g}^{\sharp} \leq c_{\Lambda,f,g} \leq C_{\Lambda,f,g} \leq C_{\Lambda,f,g}^{\flat}$.

\begin{prop}\label{positivity} Consider a non-constant symbol $f\in C^0({\mathbb T}^d)$ with $|f|=1$ and $g=(g_j)_{j=1}^d \subset C^2({\mathbb T}^d)$ a family of real-valued functions. Suppose there exists an open set $\Theta \subset {\mathbb T}^d$, which contains the support of $g$ and on which $f$ is continuously differentiable. Let $\Lambda\subset {\mathbb T}$ be a Borel set such that $e^{i\Lambda} \cap$ Ran $f \neq \emptyset$. Then,
\begin{equation*}
C_{\Lambda,f,(-i\bar{f}\nabla f\cdot g)}E_{\Lambda}(L_f) \geq E_{\Lambda}(L_f)(L_f^* A_g L_f -A_g)E_{\Lambda}(L_f) \geq c_{\Lambda,f,(-i\bar{f}\nabla f\cdot g)}E_{\Lambda}(L_f)  \enspace .
\end{equation*}
In particular, for any real-valued function $\phi\in C^0({\mathbb T})$ vanishing outside $\Lambda,$
\begin{equation*}
C_{\Lambda,f,(-i\bar{f}\nabla f\cdot g)}\Phi(L_f)^2 \geq \Phi(L_f)(L_f^* A_g L_f -A_g)\Phi(L_f) \geq c_{\Lambda,f,(-i\bar{f}\nabla f\cdot g)}\Phi(L_f)^2  \enspace .
\end{equation*}
\end{prop}
\pf Note that $E_{\Lambda}(L_f)=\chi_{f^{-1}(\Lambda)}(L_f)$ for any Borel set $\Lambda \subset {\mathbb R}$. By Lemma \ref{H0-4}, $L_f\in C^1(A_g)$ and $(L_f^* A_g L_f -A_g)=(L_{\bar{f}} A_g L_f -A_g)=L_{-i\bar{f}\nabla f \cdot g}$, hence the result. \ep

\begin{prop}\label{positivity1} Consider a non-constant symbol $f\in C^3({\mathbb T}^d)$ with $|f|=1$ and $g=(g_j)_{j=1}^d \subset C^2({\mathbb T}^d)$ a family of real-valued functions. Suppose there exists an open set $\Theta \subset {\mathbb T}^d$, which contains the support of $g$ and on which $f$ is continuously differentiable. Let $\Lambda\subset {\mathbb T}$ be a Borel set such that $e^{i\Lambda} \cap$ Ran $f \neq \emptyset$. Let $U\in C^1(A_g)$ be a unitary operator defined on $l^2({\mathbb Z}^d)$ such that $L_f^*U-I$ and $\mathrm{ad}_{A_g} L_f^*U$ are compact. Then, for any real-valued function $\phi\in C^0({\mathbb T})$ vanishing outside $\Lambda,$
\begin{equation}\label{est}
C_{\Lambda,f,(-i\bar{f}\nabla f\cdot g)}\Phi(U)^2 \gtrsim \Phi(U)(U^* A_g U -A_g)\Phi(U) \gtrsim c_{\Lambda,f,(-i\bar{f}\nabla f\cdot g)}\Phi(U)^2 \enspace .
\end{equation}
We also have that for any open set $\Lambda'$ such that $\overline{\Lambda}\subset \Lambda',$
\begin{equation}\label{mp1}
C_{\Lambda',f,(-i\bar{f}\nabla f\cdot g)}E_{\Lambda}(U) \gtrsim E_{\Lambda}(U) (U^* A_g U -A_g) E_{\Lambda}(U) \gtrsim c_{\Lambda',f,(-i\bar{f}\nabla f\cdot g)} E_{\Lambda}(U) \enspace .
\end{equation}
\end{prop}
\pf (\ref{est}) follows from the combination of Lemma \ref{prop-compact} with Proposition \ref{positivity}. Now, fix an open set $\Lambda'$ such that $\overline{\Lambda}\subset \Lambda'$. Apply (\ref{est}) to $\Lambda'$ and any real-valued function $\phi \in C^0({\mathbb T})$ vanishing on ${\mathbb T}\setminus {\Lambda'}$, such that $\phi \upharpoonright \overline{\Lambda}\equiv 1$ (Urysohn Lemma). After multiplying the corresponding inequalities on the left and right by $E_{\Lambda}(U)$, we deduce (\ref{mp1}). \ep

Let us reformulate Proposition \ref{positivity1} in our context.
\begin{cor}\label{positivity2} Consider a non-constant symbol $f\in C^3({\mathbb T}^d)$ with $|f|=1$. Let $\Lambda$ be a Borel set such that $e^{i\Lambda} \subset$ Ran $f$. Let $U\in C^1(A_{if\nabla \bar{f}})$ be a unitary operator defined on $l^2({\mathbb Z}^d)$ such that $L_f^*U-I$ and $\mathrm{ad}_{if \nabla \bar{f}} L_f^*U$ are compact. Then, for any real-valued function $\phi\in C^0({\mathbb T})$ vanishing outside $\Lambda,$
\begin{equation}\label{mp2a}
C_{\Lambda,f,|\nabla f|^2}\Phi(U)^2 \gtrsim \Phi(U)(U^* A_{if\nabla \bar{f}} U -A_{if\nabla \bar{f}})\Phi(U) \gtrsim c_{\Lambda,f,|\nabla f|^2}\Phi(U)^2  \enspace .
\end{equation}
We also have that for any open set $\Lambda'$ such that $\overline{\Lambda}\subset \Lambda',$
\begin{equation}\label{mp2b}
C_{\Lambda',f, |\nabla f|^2}E_{\Lambda}(U) \gtrsim E_{\Lambda}(U) (U^* A_{i f\nabla \bar{f}} U -A_{i f\nabla \bar{f}}) E_{\Lambda}(U) \gtrsim c_{\Lambda',f, |\nabla f|^2} E_{\Lambda}(U)  \enspace .
\end{equation}
If $e^{i\overline{\Lambda}} \subset$ Ran $f\setminus f(\kappa_f)$, then $c_{\Lambda,f,|\nabla f|^2}>0$ and there exists an open set $\Lambda'$ containing $\overline{\Lambda}$ such that $c_{\Lambda',f,|\nabla f|^2}>0$.
\end{cor}
\pf (\ref{mp2a}) and (\ref{mp2b}) follow respectively from (\ref{est}) and (\ref{mp1}) by choosing $g=-if\nabla\bar{f}$. Now, assume $e^{i\Lambda}\subset$ Ran $f\setminus f(\kappa_f)$. Since $|\nabla f|^2$ is a continuous function on ${\mathbb T}^d$ and $\overline{f^{-1}(e^{i\Lambda})}\subset f^{-1}(e^{i\overline{\Lambda}})\subset {\mathbb T}^d \setminus f^{-1}(f(\kappa_f)) \subset {\mathbb T}^d \setminus \kappa_f$, we deduce that $c_{\Lambda,f,|\nabla f|^2}>0$. Moreover, we can pick an open set $\Lambda'$ such that $\overline{\Lambda} \subset \Lambda'$ and $e^{i\overline{\Lambda'}}\subset$ Ran $f \setminus f(\kappa_f)$. According to what we have just proved, $c_{\Lambda',f,|\nabla f|^2}>0$. \ep

\begin{cor}\label{positivity3} Consider a non-constant symbol $f\in C^0({\mathbb T}^d)$ with $|f|=1$. Let $\Lambda \subset {\mathbb T}$ be a Borel set which is $M_f$-good and $\eta$ a $(M_f,\Lambda)$-adapted smooth real-valued function. Let $U\in C^1(A_{i\eta f\nabla \bar{f}})$ be a unitary operator defined on $l^2({\mathbb Z}^d)$ such that $L_f^*U-I$ and $\mathrm{ad}_{i\eta f \nabla \bar{f}} L_f^*U$ are compact. Then, for any real-valued function $\phi\in C^0({\mathbb T})$ vanishing outside $\Lambda,$
\begin{equation}\label{mp3}
C_{\Lambda,f,\eta |\nabla f|^2}\Phi(U)^2 \gtrsim \Phi(U)(U^* A_{i\eta f\nabla \bar{f}} U -A_{i\eta f\nabla \bar{f}})\Phi(U) \gtrsim c_{\Lambda,f,\eta |\nabla f|^2}\Phi(U)^2
\end{equation}
where $c_{\Lambda,f,\eta |\nabla f|^2}>0$.
\end{cor}
\pf (\ref{mp3}) follows from (\ref{est}) by choosing $g=-i\eta f\nabla\bar{f}$. The function $\eta |\nabla f|^2$ is continuous on ${\mathbb T}^d$. Since $e^{i\overline{\Lambda}}\subset$ Ran $f\setminus f(\kappa_f)$, then $\overline{f^{-1}(e^{i\Lambda})}\subset f^{-1}(e^{i\overline{\Lambda}})\subset {\mathbb T}^d \setminus f^{-1}(f(\kappa_f)) \subset {\mathbb T}^d \setminus \kappa_f$ and $c_{\Lambda,f,\eta |\nabla f|^2}>0$. \ep

\subsection{Proof of Theorem \ref{laurent-glob}}

By Lemma \ref{H0-4}, $L_f\in C^3(A_{if\nabla \bar{f}})\subset {\cal C}^{1,1}(A_{if\nabla \bar{f}})$. Since $U\in {\cal C}^{1,1}(A_{if\nabla \bar{f}})$, $(L_f^*U-I) \in {\cal C}^{1,1}(A_{if\nabla \bar{f}})$. Because $(L_f^*U-I)$ is compact, $\mathrm{ad}_{A_{if\nabla \bar{f}}}L_f^*U= \mathrm{ad}_{A_{if\nabla \bar{f}}}(L_f^*U-I)$ is also compact by Lemma \ref{ad-compact}. Due to Corollary \ref{positivity2}, $U$ is propagating w.r.t. $A_{if\nabla \bar{f}}$ on any Borel set $\Lambda$ such that $e^{i\overline{\Lambda}}\subset$ Ran $f\setminus f(\kappa_f)$. Statement (a) follows from Corollary \ref{virial-2}. By requiring $\Lambda$ to be an open set in this discussion, we deduce (b) from Theorem \ref{nosc2}.

\subsection{Proof of Theorem \ref{laurent-loc}}

By Lemma \ref{H0-4}, $L_f\in C^3(A_{i\eta f\nabla \bar{f}})\subset {\cal C}^{1,1}(A_{i\eta f\nabla \bar{f}})$. Since $U\in {\cal C}^{1,1}(A_{i\eta f\nabla \bar{f}})$, $(L_f^*U-I) \in {\cal C}^{1,1}(A_{i\eta f\nabla \bar{f}})$. Because $(L_f^*U-I)$ is compact, $\mathrm{ad}_{A_{i\eta f\nabla \bar{f}}}L_f^*U$ is also compact by Lemma \ref{ad-compact}. Let $\Lambda'\subset {\mathbb T}$ be a Borel set such that $\overline{\Lambda'}\subset \Lambda$. Apply (\ref{mp3}) to any real-valued function $\phi\in C^0({\mathbb T})$, vanishing outside $\Lambda$ such that $\phi \upharpoonright \overline{\Lambda'}\equiv 1$. Multiplying the corresponding inequalities on the left and right by $E_{\Lambda'}(U)$ entails
\begin{equation*}
C_{\Lambda,f, \eta |\nabla f|^2}E_{\Lambda'}(U) \gtrsim E_{\Lambda'}(U) (U^* A_{i \eta f\nabla \bar{f}} U -A_{i \eta f\nabla \bar{f}}) E_{\Lambda'}(U) \gtrsim c_{\Lambda,f, \eta |\nabla f|^2} E_{\Lambda'}(U)
\end{equation*}
where $c_{\Lambda,f, \eta |\nabla f|^2}>0$. Statement (a) follows from Corollary \ref{virial-2}. Now, by requiring $\Lambda'$ to be an open set in this discussion, we deduce (b) from Theorem \ref{nosc2}.

\subsection{Proof of Proposition \ref{transl}}

We note that the symbol $f$ is analytic and has no critical point ($\kappa_f=\emptyset$). According to Lemma \ref{H0-4}, $L_f=T^{\alpha}$ belongs to $C^{\infty}(\sum_{j=1}^d \alpha_j X_j)$ and $(T^{-\alpha} (\sum_{j=1}^d \alpha_j X_j)T^{\alpha} -(\sum_{j=1}^d \alpha_j X_j))=|\alpha|^2$ ($|\alpha |^2=\sum_{j=1}^d \alpha_j^2$). Since $U\in C^1(\sum_{j=1}^d \alpha_j X_j)$ (recall that ${\cal C}^{1,1}(\sum_{j=1}^d \alpha_j X_j)\subset C^1(\sum_{j=1}^d \alpha_j X_j)$) and $\mathrm{ad}_{\sum_{j=1}^d \alpha_j X_j}T^{-\alpha}U$ is compact we deduce from Lemma \ref{prop-compact} that:
\begin{equation}\label{mourre-tr}
\left(U^* (\sum_{j=1}^d \alpha_j X_j)U-(\sum_{j=1}^d \alpha_j X_j)\right)\simeq |\alpha|^2 \enspace.
\end{equation}
Since $\alpha\neq 0$, the conclusions are deduced by applying Corollary \ref{virial-2} and Theorem \ref{nosc2}.

\section{Propagation estimates}

The commutator formalism also allows to derive some propagation estimates. In this section, $f$ belongs to $C^3({\mathbb T})$, $|f|=1$. We will also use freely the notations introduced in Section 4 and write for any $\varphi \in {\cal D}_X,$
\begin{equation}\label{normX}
\|\varphi\|_X = \sqrt{\|\varphi\|^2+\sum_{j=1}^d\|X_j \varphi \|^2} \enspace .
\end{equation}
We consider first the propagation properties generated by the iterations of the unitary operator $L_f$:
\begin{prop}\label{laurent-1} Consider a non-constant symbol $f\in C^3({\mathbb T}^d)$ with $|f|=1$. Then, for any $(\varphi,\psi) \in {\cal H}\times {\cal D}(A_{if\nabla\bar{f}}),$
\begin{equation}\label{propagfree1}
\lim_{n\rightarrow \pm \infty} n^{-1} \bra L_f^n\varphi, A_{if\nabla\bar{f}} L_f^n\psi \ket = \bra \varphi, L_{|\nabla f|^2}\psi \ket \geq 0 \enspace ,
\end{equation}
and for any $\psi \in {\cal D}_X,$
\begin{equation}\label{propagfree2}
\lim_{n\rightarrow \pm\infty} |n|^{-1} \|L_f^n\psi\|_X = \sqrt{\bra \psi, L_{|\nabla f|^2}\psi \ket} = \|L_{|\nabla f|}\psi \| \enspace .
\end{equation}
\end{prop}
For a similar result in the Hamiltonian case, see e.g. \cite{ak}. These estimates are preserved in a weaker form when considering suitable unitary perturbations $U$ of $L_f$. We have that:
\begin{prop}\label{laurent-2ub} Consider a non-constant symbol $f\in C^3({\mathbb T}^d)$ with $|f|=1$. Let $U\in \cap_{j=1}^d C^1(X_j) \cap C^1 (A_{if\nabla\bar{f}})$ be a unitary operator defined on $l^2({\mathbb Z}^d)$ such that $L_f^*U \in \cap_{j=1}^d C^1(X_j^2)$. Then, for any $\psi \in {\cal D}_X,$
\begin{equation}\label{limsup(b)}
\limsup_{n\rightarrow \pm\infty} |n|^{-1} \|U^n\psi\|_X \leq \sqrt{\|\mathrm{ad}_{A_{if\nabla \bar{f}}}U \|} \|\psi\| \enspace .
\end{equation}
\end{prop}
See Section 5.3 for the proof.

\begin{prop}\label{laurent-2lb} Consider a non-constant symbol $f\in C^3({\mathbb T}^d)$ with $|f|=1$. Let $U\in \cap_{j=1}^d C^1(X_j) \cap C^1 (A_{if\nabla\bar{f}})$ be a unitary operator defined on $l^2({\mathbb Z}^d)$ such that $L_f^*U \in \cap_{j=1}^d C^1(X_j^2)$. Assume that $\mathrm{ad}_{A_{if\nabla\bar{f}}} L_f^*U$ and $L_f^*U-I$ are compact. Let $\Lambda \subset {\mathbb T}$ be a Borel set such that $e^{i\Lambda}\subset$ Ran $f$.
\begin{itemize}
\item For any $\psi \in {\cal H}$, such that $E_{\Lambda}(U)\psi \in {\cal H}_c(U)\cap {\cal D}_X,$
\begin{eqnarray}
\sqrt{c_{\Lambda,f,|\nabla f|^2}^{\sharp}} \|E_{\Lambda}(U) \psi\| &\leq & \liminf_{n\rightarrow \pm\infty} |n|^{-1} \|U^n E_{\Lambda}(U) \psi\|_X \nonumber \\
\limsup_{n\rightarrow \pm\infty} |n|^{-1} \|U^n E_{\Lambda}(U)\psi\|_X &\leq & \sqrt{C_{\Lambda,f,|\nabla f|^2}^{\flat}} \|E_{\Lambda}(U) \psi\| \enspace . \label{liminf(b)}
\end{eqnarray}
\item Let $\phi \in C^1({\mathbb T})$ be real-valued, vanishing outside $\Lambda$ with $\phi' \in {\cal A}({\mathbb T})$. For any $\psi \in {\cal H}_c(U)\cap {\cal D}_X,$
\begin{eqnarray}
\sqrt{c_{\Lambda,f,|\nabla f|^2}} \|\Phi(U) \psi\| &\leq & \liminf_{n\rightarrow \pm\infty} |n|^{-1} \|U^n \Phi(U) \psi\|_X \nonumber \\
\limsup_{n\rightarrow \pm\infty} |n|^{-1} \|U^n \Phi(U) \psi\|_X & \leq & \sqrt{C_{\Lambda,f,|\nabla f|^2}} \|\Phi(U) \psi\| \enspace . \label{liminf(c)}
\end{eqnarray}
\item If $e^{i\overline{\Lambda}} \subset$ Ran $f\setminus f(\kappa_f)$, then $c_{\Lambda,f,|\nabla f|^2}\geq c_{\Lambda,f,|\nabla f|^2}^{\sharp}>0$.
\end{itemize}
\end{prop}
See Section 5.4 for the proof. We note that if $U\in {\cal C}^{1,1}(A_{if\nabla\bar{f}})$ and if $L_f^*U$ is compact, then $\mathrm{ad}_{A_{if\nabla\bar{f}}} L_f^*U$ is compact (see Lemma \ref{ad-compact}).

\noindent{\bf Remark:} The construction of a non trivial vector $\varphi$ satisfying the condition $E_{\Lambda}(U)\varphi \in {\cal H}_c(U)\cap {\cal D}_X$ in Proposition \ref{laurent-2lb}, can be easily performed in the following situation. Let $\Lambda\subset {\mathbb T}$ be an open interval such that $e^{i\overline{\Lambda}} \subset$ Ran $f\setminus f(\kappa_f)$. If $U\in C^1(A_{if\nabla \bar{f}})$, it follows from Corollaries \ref{virial-2} and \ref{positivity2} that $\sigma_{pp}(U)\cap e^{i\Lambda}$ is finite. Moreover $U\in \cap_{j=1}^d C^1(X_j)$, so $\Phi(U)$ belong to $\cap_{j=1}^d C^1(X_j)$ for any $\phi \in C^1({\mathbb T})$ with $\phi' \in {\cal A}({\mathbb T})$ (see e.g. Proposition \ref{10}). By considering such a function $\phi$, real-valued and vanishing outside $\Lambda'$ where $e^{i\Lambda'}= e^{i\Lambda}\setminus \sigma_{pp}(U)$, we define $\varphi=\Phi(U)\psi$ with $\psi \in {\cal D}_X$ and have that $E_{\Lambda}(U)\varphi=\varphi \in {\cal D}_X \cap {\cal H}_c(U)$.

In the translation case, Proposition \ref{laurent-2lb} can be simplified as follows:
\begin{prop}\label{laurent-2tr} Let $\alpha\in {\mathbb Z}^d\setminus \{0\}$. Let $U\in \cap_{j=1}^d C^1(X_j) \cap C^1(\sum_{j=1}^d \alpha_j X_j)$ be a unitary operator defined on $l^2({\mathbb Z}^d)$ such that $T^{-\alpha}U \in \cap_{j=1}^d C^1(X_j^2)$. Assume that $\mathrm{ad}_{\sum_{j=1}^d \alpha_j X_j} T^{-\alpha}U$ is compact. Then for any $\psi \in {\cal H}_c(U) \cap {\cal D}_X,$
\begin{equation}\label{lim(b)}
\lim_{n\rightarrow \pm \infty} |n|^{-1} \| U^n \psi\|_X = |\alpha| \|\psi \| \enspace .
\end{equation}
\end{prop}
See Section 5.5 for the proof.

Sections 5.1 and 5.2 allow to articulate the proofs of Propositions \ref{laurent-2ub}, \ref{laurent-2lb} and \ref{laurent-2tr} with Mourre Theory.

\subsection{Preliminaries}

In this section, we build on the contents of Section 4.1. ${\cal H}$ denotes a fixed infinite dimensional Hilbert space while $A$ is a fixed self-adjoint operator (densely) defined on ${\cal H}$ with domain ${\cal D}(A)$. 

If a unitary operator $U$ defined on ${\cal H}$ belongs to $C^1(A)$, $U^n \in C^1(A)$ for any $n\in {\mathbb Z}$. So for any $\psi \in {\cal D}(A)$ and any $n\in {\mathbb Z}$, $U^n \psi  \in {\cal D}(A)$. Thus, according to Proposition \ref{propag07}, we have that for any $(\varphi,\psi) \in {\cal H}\times {\cal D}(A)$ and any $n\in {\mathbb N},$
\begin{eqnarray}
\langle U^n\varphi,AU^n\psi\rangle\ - \langle\varphi,A\psi\rangle &=& \sum_{m=0}^{n-1} \bra U^m\varphi, (U^*AU-A) U^m \psi \ket \label{identity23}\\
\langle U^{*n}\varphi,AU^{*n}\psi\rangle\ - \langle\varphi,A\psi\rangle &=& - \sum_{m=1}^{n} \bra U^{*m}\varphi, (U^*AU-A) U^{*m} \psi \ket \label{identity23-}
\end{eqnarray}
We recall that $(U^*AU-A)$ is the (unique) bounded operator associated to the sesquilinear form $F_+^o$ defined in Lemma \ref{propag0} and that $(U^*AU-A)=U^* \mathrm{ad}_A U$. Similar formulas can be drawn where the RHS of (\ref{identity23}) and (\ref{identity23-}) are expressed in terms of $(A-UAU^*)$ instead of $(U^*AU-A)$. It follows that:

\begin{prop}\label{laurent-0} Let $U$ be a unitary operator in $C^1(A)$. Then, for any $(\varphi,\psi) \in {\cal H}\times {\cal D}(A),$
\begin{equation*}
\limsup_{n\rightarrow \pm \infty} |n|^{-1} |\bra U^n\varphi, A U^n\psi \ket | \leq C \|\varphi\| \|\psi\| \enspace ,
\end{equation*}
where $C=\|(U^*AU-A)\|=\|\mathrm{ad}_A U\|=\|\mathrm{ad}_A U^*\|=\|(A-UAU^*)\|$.
\end{prop}
\noindent{\bf Proof:} According to Lemma \ref{propag01} (recall that $\|U\|=\|U^{-1}\|=1$), if $U$ belongs to $C^1(A)$, then for all $n\in {\mathbb N}$, $| \bra U^n\varphi,AU^n \psi \ket - \bra \varphi, A\psi \ket | \leq C n \|\varphi\| \|\psi\|$ and $| \bra \varphi, A\psi \ket - \bra U^{*n}\varphi,AU^{*n} \psi \ket | \leq C n \|\varphi\| \|\psi\|$, which proves the result. \ep

\noindent{\bf Remark:} If a unitary operator $U$ belongs to $C^1(A)$ and commutes with $(U^*AU-A)$, then
\begin{equation}\label{identitysum1}
\langle U^n\varphi,AU^n\psi\rangle\ = \langle\varphi,A\psi\rangle + n \bra \varphi, (U^*AU-A) \psi \ket
\end{equation}
for any $(\varphi,\psi) \in {\cal H}\times {\cal D}(A)$ and for all $n\in {\mathbb Z}$.

\noindent{\bf Remark:} Let $U\in C^1(A)$. We claim that $U$ commutes with $(U^*AU-A)$ if and only if it commutes with $\mathrm{ad}_A U$. Note that given any $B\in {\cal B}({\cal H})$: $0=[B,U^{-1}U]=U^{-1}[B,U]+[B,U^{-1}]U$. This means that if the unitary operator $U$ belongs to $C^1(A)$, $[(U^*AU-A),U]=0$ if and only if $[(U^*AU-A),U^*]=0$. On the other hand, $U$ is normal and $(U^*AU-A)=U^*(\mathrm{ad}_A U)$ (see e.g. Lemma \ref{propag0}), which implies our claim. 

Let $U$ be a unitary operator acting on ${\cal H}$ and $\varphi\in {\cal H}_c(U)$ where ${\cal H}_c(U)$ denotes the continuous subspace of $U$. It follows from RAGE Theorem (see e.g. Theorem 3.2 in \cite{ev}) that:
\begin{equation}\label{rage}
\lim_{n\rightarrow \infty} \frac{1}{n} \sum_{m=0}^{n-1}\|KU^m\varphi\|=\lim_{n\rightarrow \infty} \frac{1}{n} \sum_{m=1}^{n}\|KU^{*m}\varphi\|=0\enspace .
\end{equation}
for any compact operator $K$. This allows us to prove that:
\begin{lem}\label{pp2c0} Let $U$ be a unitary operator which belongs to $C^1(A)$ such that: $(U^*AU-A)=B+K$ where $B$ is a bounded operator that commutes with $U$ and $K$ is compact. Let $(\varphi,\psi)\in {\cal H}\times {\cal D}(A)$ and assume that either $\varphi$ or $\psi$ belong to ${\cal H}_c(U)$. Then,
\begin{equation*}
\lim_{n\rightarrow \pm \infty} n^{-1} \bra U^n \varphi, A U^n \psi \ket = \bra \varphi, B\psi \ket \enspace .
\end{equation*}
\end{lem}
\noindent{\bf Proof:} First note that (\ref{identity23}) and (\ref{identity23-}) rewrite: for all $n\in {\mathbb N},$
\begin{eqnarray*}
\langle U^n\varphi,AU^n\psi\rangle\ - \langle\varphi,A\psi\rangle &=& n \bra \varphi, B\psi\ket +\sum_{m=0}^{n-1} \bra U^m\varphi, K U^m \psi \ket \\
\langle U^{*n}\varphi,AU^{*n}\psi\rangle\ - \langle\varphi,A\psi\rangle &=& -n \bra \varphi, B\psi\ket - \sum_{m=1}^{n} \bra U^{*m}\varphi, K U^{*m} \psi \ket \enspace .
\end{eqnarray*}
The conclusion follows from (\ref{rage}) since $K$ and $K^*$ are simultaneously compact. \ep

In a more general setup, let $U\in C^1(A)$ and let us assume there exist two operators $B$ and $K$, respectively bounded and compact such that $(U^*AU-A)\gtreqless B+K$. Since for all $n\in {\mathbb Z}$, $U^n$ leaves ${\cal H}_c(U)\cap {\cal D}(A)$ invariant, it follows from (\ref{identity23}), (\ref{identity23-}) and (\ref{rage}) that for all $\varphi \in {\cal H}_c(U)\cap {\cal D}(A),$
\begin{eqnarray*}
\liminf_{n\rightarrow \infty} n^{-1} \bra U^n\varphi, A U^n\varphi \ket &\gtreqless & \liminf_{n\rightarrow \infty} n^{-1} \sum_{m=0}^{n-1}\bra U^m\varphi, B U^m\varphi \ket \\
\liminf_{n\rightarrow -\infty} n^{-1} \bra U^n\varphi, A U^n\varphi \ket &\gtreqless & \liminf_{n\rightarrow -\infty} |n|^{-1} \sum_{m=1}^{|n|}\bra U^{*m}\varphi, B U^{*m}\varphi \ket\\
\limsup_{n\rightarrow \infty} n^{-1} \bra U^n\varphi, A U^n\varphi \ket &\gtreqless & \limsup_{n\rightarrow \infty} n^{-1} \sum_{m=0}^{n-1}\bra U^m\varphi, B U^m\varphi \ket \\
\limsup_{n\rightarrow -\infty} n^{-1} \bra U^n\varphi, A U^n\varphi \ket &\gtreqless & \limsup_{n\rightarrow -\infty} |n|^{-1} \sum_{m=1}^{|n|}\bra U^{*m}\varphi, B U^{*m}\varphi \ket  \enspace .
\end{eqnarray*}
Following the same idea, we have also that:
\begin{prop}\label{pp2c} Let $\Theta\subset {\mathbb T}$ be a Borel set, $\phi \in C^0({\mathbb T})$ and $U\in C^1(A)$.
\begin{itemize}
\item If $c E_{\Theta}(U) \lesssim E_{\Theta}(U) (U^*AU-A) E_{\Theta}(U)$ for some $c\in {\mathbb R}$, then for any $\varphi \in {\cal H}$ such that $E_{\Theta}(U)\varphi \in {\cal H}_c(U)\cap {\cal D}(A),$
\begin{equation*}
c \|E_{\Theta}(U)\varphi \|^2 \leq \liminf_{n\rightarrow \pm \infty} n^{-1} \bra U^n E_{\Theta}(U) \varphi, A U^n E_{\Theta}(U) \varphi \ket \enspace .
\end{equation*}
\item If $E_{\Theta}(U) (U^*AU-A) E_{\Theta}(U)\lesssim C E_{\Theta}(U)$ for some $C\in {\mathbb R}$, then for any $\varphi \in {\cal H}$ such that $E_{\Theta}(U)\varphi \in {\cal H}_c(U)\cap {\cal D}(A),$
\begin{equation*}
\limsup_{n\rightarrow \pm \infty} n^{-1} \bra U^n E_{\Theta}(U) \varphi, A U^n E_{\Theta}(U) \varphi \ket \leq C \|E_{\Theta}(U)\varphi\|^2 \enspace .
\end{equation*}
\item If $c \Phi(U)^2 \lesssim \Phi(U) (U^*AU-A) \Phi(U)$ for some $c\in {\mathbb R}$, then for any $\varphi \in {\cal H}$ such that $\Phi(U)\varphi \in {\cal H}_c(U)\cap {\cal D}(A),$
\begin{equation*}
c \|\Phi(U)\varphi \|^2 \leq \liminf_{n\rightarrow \pm \infty} n^{-1} \bra U^n \Phi(U) \varphi, A U^n \Phi(U) \varphi \ket \enspace .
\end{equation*}
\item If $\Phi(U) (U^*AU-A) \Phi(U)\lesssim C \Phi(U)^2$ for some $C\in {\mathbb R}$, then for any $\varphi \in {\cal H}$ such that $\Phi(U)\varphi \in {\cal H}_c(U)\cap {\cal D}(A),$
\begin{equation*}
\limsup_{n\rightarrow \pm \infty} n^{-1} \bra U^n \Phi(U) \varphi, A U^n \Phi(U) \varphi \ket \leq C \|\Phi(U)\varphi\|^2 \enspace .
\end{equation*}
\end{itemize}
\end{prop}

\subsection{Preliminaries. Continued}

We come back to the context and notations of Section 2. 
\begin{prop}\label{propag6} Let $U$ be a unitary operator defined on $l^2({\mathbb Z}^d)$. Assume that $U$ belongs to $\cap_{j=1}^d C^1(X_j)$ and write: $C_j:=(X_j-UX_jU^*)$, $j\in \{1, \ldots , d\}$. Then, for any $\psi \in {\cal D}_X$ and all $n\in {\mathbb N},$
\begin{eqnarray}
\|U^n \psi\|_X^2 - \|\psi\|_X^2 &=& 2 \sum_{m=1}^n \sum_{j=1}^d \Re \bra C_j U^m\psi, X_j U^m\psi\ket - \sum_{m=1}^n \sum_{j=1}^d \| C_j U^m\psi\|^2 \label{xnorma-sum}\\
\|U^{*n} \psi\|_X^2 - \|\psi\|_X^2 &=& - 2 \sum_{m=0}^{n-1} \sum_{j=1}^d \Re \bra C_j U^{*m}\psi, X_j U^{*m}\psi\ket + \sum_{m=0}^{n-1} \sum_{j=1}^d \| C_j U^{*m}\psi\|^2 \enspace . \label{xnorma-sum2}
\end{eqnarray}
If $\limsup_{n\rightarrow \pm \infty} |n|^{-1} |\sum_{j=1}^d \Re \bra C_j U^n\psi,X_j U^n\psi\ket |< \infty$, then
$$
\limsup_{n\rightarrow \pm \infty} |n|^{-1} \|U^n \psi\|_X \leq \left(\limsup_{n\rightarrow \pm \infty} |n|^{-1} |\sum_{j=1}^d \Re \bra C_j U^n\psi,X_j U^n\psi\ket |\right)^{\frac{1}{2}} \enspace .
$$
If $0\leq \liminf_{n\rightarrow \pm \infty} n^{-1} \sum_{j=1}^d\Re \bra C_j U^n\psi, X_j U^n\psi\ket \leq \limsup_{n\rightarrow \pm \infty} n^{-1} \sum_{j=1}^d\Re \bra C_j U^n\psi, X_j U^n\psi\ket < \infty$, then
\begin{eqnarray*}
\sqrt{\liminf_{n\rightarrow \pm \infty} n^{-1} \sum_{j=1}^d \Re \bra C_j U^n\psi,X_j U^n\psi\ket} &\leq & \liminf_{n\rightarrow \pm \infty} |n|^{-1} \|U^n \psi\|_X \\
\limsup_{n\rightarrow \pm \infty} |n|^{-1} \|U^n \psi\|_X &\geq & \sqrt{\limsup_{n\rightarrow \pm \infty} n^{-1} \sum_{j=1}^d \Re \bra C_j U^n\psi,X_j U^n\psi\ket} \enspace .
\end{eqnarray*}
\end{prop}
\noindent{\bf Proof:} If $U\in \cap_{j=1}^d C^1(X_j)$ then for any vector $\psi \in {\cal D}_X$, the vectors $U^n\psi$ also belong to ${\cal D}_X$ for all $n\in {\mathbb Z}$. The operators $C_j$ are bounded for all $j\in \{1,\ldots, d\}$.  Now, fix $j\in \{1,\ldots, d\}$ and $\psi\in {\cal D}_X$. One has that for all $n\in {\mathbb N},$
\begin{eqnarray*}
\|X_j U^n \psi\|^2 - \|X_j\psi\|^2 &=& \sum_{m=1}^n \| X_j U^m\psi\|^2 - \|UX_jU^* U^{m}\psi\|^2 \\
&=& 2 \sum_{m=1}^n \Re \bra C_j U^m\psi, X_j U^m\psi\ket - \sum_{m=1}^n \|C_j U^m\psi\|^2 \\
\|X_j U^{*n} \psi\|^2 - \|X_j\psi\|^2 &=& \sum_{m=0}^{n-1} \| UX_j U^{*} U^{*m}\psi\|^2 - \|X_j U^{*m}\psi\|^2 \\
&=&- 2 \sum_{m=0}^{n-1} \Re \bra C_j U^{*m}\psi, X_j U^{*m}\psi\ket + \sum_{m=0}^{n-1} \| C_j U^{*m}\psi\|^2 \enspace .
\end{eqnarray*}
Note that $\sup_{m\in {\mathbb Z}}\| C_j U^m\psi\|^2 \leq \|C_j\|^2 \|\psi \|^2 $. The conclusions follow after summing over $j\in \{1,\ldots,d\}$ in the previous identities. \ep

\begin{cor}\label{propag7} Let $U$ be a unitary operator defined on $l^2({\mathbb Z}^d)$. Assume that $U$ belongs to $\cap_{j=1}^d C^1(X_j)$ and that $U$ commutes with $C_j:=(X_j-UX_jU^*)$ (or equivalently with $\mathrm{ad}_{X_j} U$) for all $j\in \{1,\ldots,d\}$. Then, for any $\psi \in {\cal D}_X$ and all $n\in {\mathbb Z},$
\begin{equation}\label{normX2}
\|U^n \psi\|_X^2 - \|\psi\|_X^2 = n^2 \sum_{j=1}^d \|C_j \psi\|^2 + 2n \sum_{j=1}^d \Re \bra C_j \psi, X_j \psi\ket \enspace .
\end{equation}
\end{cor}
\noindent{\bf Proof:} For all $j\in \{1,\ldots,d\}$, all $\psi \in {\cal D}_X$ and all $m\in {\mathbb Z}$, we have that: $\|C_jU^m\psi\|=\|C_j\psi\|$ and
\begin{eqnarray*}
\bra C_j U^m \psi, X_jU^m\psi \ket &=& \bra U^m C_j \psi,X_jU^m\psi \ket = \bra C_j \psi, X_j\psi \ket + m \bra  C_j \psi, U^*(X_j-UX_jU^*)U\psi \ket  \\
&= &  \bra C_j \psi, X_j\psi \ket + m \|C_j \psi\|^2 \enspace .
\end{eqnarray*}
The result follows by plugging the former identity into (\ref{xnorma-sum}) and (\ref{xnorma-sum2}). \ep

The proofs of Propositions \ref{laurent-2ub} and \ref{laurent-2lb} require also the following result:
\begin{lem}\label{qf} Let ${\cal H}$ be a Hilbert space and $A$ a self-adjoint operator with domain ${\cal D}(A)\subset {\cal H}$. Let $U$ be a unitary operator defined on ${\cal H}$. If $U$ belongs to $C^1(A)\cap C^1(A^2)$, then the sesquilinear forms $Q_+$ and $Q_-$ defined on ${\cal D}(A)\times {\cal D}(A)$ by
\begin{eqnarray*}
Q_+(\varphi,\psi) &=& \bra (U^*AU-A) \varphi, A\psi\ket + \bra A\varphi, (U^*AU-A) \psi \ket \\
Q_-(\varphi,\psi) &=& \bra (A-UAU^*) \varphi, A\psi\ket + \bra A\varphi, (A-UAU^*) \psi \ket
\end{eqnarray*}
extend as bounded forms on ${\cal H}\times {\cal H}$.
\end{lem}
\noindent{\bf Proof:} First, we observe that ${\cal D}(A^2)$ is a core for $A$ (see \cite{schmu} for details). Second, $U{\cal D}(A)={\cal D}(A)$ and $U{\cal D}(A^2)={\cal D}(A^2)$. So, it is enough to see that for all $(\varphi, \psi)\in {\cal D}(A^2)\times {\cal D}(A^2)$
\begin{eqnarray*}
Q_+(\varphi,\psi) &=& \bra \varphi, (U^* A^2 U-A^2)\psi \ket - \bra [A,U]\varphi, [A,U]\psi \ket \\
&=& \bra \varphi, (U^* A^2 U-A^2)\psi \ket - \bra (U^*AU-A)\varphi, (U^*AU-A)\psi \ket \\
Q_-(\varphi,\psi) &=& \bra \varphi, (A^2-UA^2U^*)\psi \ket - \bra [A,U^*]\varphi, [A,U^*]\psi \ket \\
&=& \bra \varphi, (A^2-UA^2U^*)\psi \ket - \bra (A-UAU^*)\varphi, (A-UAU^*)\psi \ket \enspace .
\end{eqnarray*}
\ep

This allows us to prove that:
\begin{lem}\label{shortcut} Let $f\in C^1({\mathbb T})$, $|f|=1$. Let $U\in \cap_{j=1}^d C^1(X_j)$ be a unitary operator defined on $l^2({\mathbb Z}^d)$ such that $W:=L_f^*U \in \cap_{j=1}^d C^1(X_j^2)$. Then, $W \in \cap_{j=1}^d C^1(X_j)$ and the sesquilinear form $Q_0$ defined on ${\cal D}_X \times {\cal D}_X$ by:
\begin{equation*}
Q_0(\varphi, \psi) := \sum_{j=1}^d \bra (X_j-WX_jW^*)L_f^*\varphi,L_f^*X_j \psi \ket 
\end{equation*}
is continuous for the topology induced by ${\cal H} \times {\cal H}$.
\end{lem}
\pf It follows from Lemma \ref{H0-4}, that $L_f$ (and $L_f^*$) belongs to $\cap_{j=1}^d C^1(X_j)$. In particular, $L_f^* {\cal D}_X \subset {\cal D}_X$ (Proposition \ref{2}) and $W\in \cap_{j=1}^d C^1(X_j)$ (Proposition \ref{3}). This allows us to define for $j\in \{1,\ldots,d\}$ the bounded operators $C_{j,W}:= (X_j-WX_jW^*)$ according to Lemma \ref{propag0} and we have that for all $(\varphi,\psi)\in {\cal D}_X \times {\cal D}_X,$
\begin{equation*}
Q_0(\varphi, \psi) = \sum_{j=1}^d \bra C_{j,W}L_f^* \psi, X_j L_f^*\psi \ket - \sum_{j=1}^d \bra C_{j,W}L_f^* \psi, [X_j, L_f^*]\psi \ket \enspace ,
\end{equation*}
where $[X_j, L_f^*] = [X_j, L_{\bar{f}}]= L_{-i\partial \bar{f}}$ for $j\in \{1,\ldots,d\}$. The fact that $W$ and $L_f$ belong to $\cap_{j=1}^d C^1(X_j)$ implies that the modulus of the last sum on the RHS is bounded by $\|\varphi \| \|\psi \|$ up to some multiplicative constant which depends on $f$ and $W$. Since $W\in \cap_{j=1}^d C^1(X_j^2)$, the modulus of the first sum on the RHS is also bounded by $\|\varphi \| \|\psi \|$ up to some multiplicative constant which depends on $W$ by Lemma \ref{qf}. This proves the lemma. \ep

\subsection{Proof of Proposition \ref{laurent-1}}

For $f\in C^3({\mathbb T}^d)$, it follows from Lemma \ref{H0-4} that $L_f\in C^1(A_{if \cdot \nabla  \bar{f}})$ and $(L_f^*A_{if \cdot \nabla  \bar{f}} L_f- A_{if \cdot \nabla  \bar{f}}) =L_{|\nabla f|^2}$. $L_{|\nabla f|^2}$ commutes with $L_f$. By (\ref{identitysum1}), we have that for all $(\varphi,\psi) \in {\cal H}\times {\cal D}(A),$
$$
\langle L_f^n\varphi,A_{if \cdot \nabla  \bar{f}} L_f^n\psi\rangle\ = \langle\varphi, A_{if \cdot \nabla  \bar{f}}\psi\rangle + n \bra\varphi, L_{|\nabla f|^2} \psi \ket
$$
which proves (\ref{propagfree1}). As shown in Section 2.2. we also have that $L_f \in \cap_{j=1}^d C^1(X_j)$ and $C_j:=(X_j-L_f X_jL_f^*)=-L_f(\mathrm{ad}_{X_j}L_f^*)=iL_{f}L_{\partial_{j}\bar{f}}$ for all $j\in \{1,\ldots,d\}$. So, $C_j$ commutes with $L_f$ for $j\in \{1,\ldots,d\}$. Since $|f|=1$, we have that for any $\psi \in {\cal H},$
$$\sum_{j=1}^d \|C_j \psi\|^2= \sum_{j=1}^d \|iL_{f}L_{\partial_{j}\bar{f}} \psi\|^2 = \sum_{j=1}^d \bra \psi, L_{|\partial_{j}f|^2} \psi\ket = \bra \psi, L_{|\nabla f|^2} \psi\ket \enspace .$$
Statement (\ref{propagfree2}) follows from (\ref{normX2}).

\subsection{Proof of Proposition \ref{laurent-2ub}}

Let us adopt some local notations: if $V$ is a unitary operator defined on $l^2({\mathbb Z}^d)$, which belongs to $C^1(X_j)$ for some $j\in \{1,\ldots,d\}$, we write $C_{j,V}=(X_j-VX_jV^*)$, according to Lemma \ref{propag0}. 

Since $U\in C^1(A_{if\nabla\bar{f}})$, we deduce from Proposition \ref{laurent-0} that for any $(\varphi,\psi) \in {\cal H}\times {\cal D}(A_{if\nabla\bar{f}}),$
\begin{equation} \label{limsup(a)}
\limsup_{n\rightarrow \pm \infty} |n|^{-1} |\bra U^n\varphi, A_{if\nabla\bar{f}} U^n\psi \ket | \leq \|\mathrm{ad}_{A_{if \cdot \nabla  \bar{f}}}U\| \|\varphi \| \|\psi \| \enspace .
\end{equation}
By Lemma \ref{H0-4} (resp. by hypothesis), $L_f$ (resp. $U$) belongs to $\cap_{j=1}^d C^1(X_j)$. So, $W:=L_f^*U$ and $U^n$ also belong to $\cap_{j=1}^d C^1(X_j)$ for any $n\in {\mathbb Z}$. For $j\in \{1, \ldots, d\}$, we can split $C_{j,U}$ as follows: $C_{j,U}= C_{j,L_{f}} +L_{f} C_{j,W}L_{f}^*$, where $C_{j,L_{f}}=L_{if\partial_j \bar{f}}$. For any $n\in {\mathbb Z}$ and any $\psi \in {\cal D}_X$, one has that $U^n\psi \in {\cal D}_X$ and
\begin{equation*}
\sum_{j=1}^d \Re \bra C_{j,U}U^n\psi, X_j U^n\psi\ket - \bra U^n\psi, A_{if\nabla\bar{f}} U^n\psi \ket = \Re \left( Q_0(U^n\psi, U^n\psi) \right)\enspace .
\end{equation*}
According to Lemma \ref{shortcut}, there exists $C>0$, such that for all $\psi \in {\cal D}_X$ and all $n\in {\mathbb Z}$, $|Q_0(U^n\psi, U^n\psi)|\leq C \|\psi\|^2$. It follows from (\ref{limsup(a)}) that for any $\psi \in {\cal D}_X$ (${\cal D}_X\subset {\cal D}(A_{if\nabla\bar{f}})$),
\begin{equation*}
\limsup_{n\rightarrow \pm \infty} |n|^{-1} \left| \sum_{j=1}^d \Re \bra C_{j,U} U^n\psi, X_j U^n\psi\ket \right| \leq \|\mathrm{ad}_{A_{if \cdot \nabla  \bar{f}}}U\| \|\psi \|^2 \enspace .
\end{equation*}
This concludes the proof in view of Proposition \ref{propag6}.

\subsection{Proof of Proposition \ref{laurent-2lb}}

Since $U\in C^1(A_{if\nabla\bar{f}})$, we deduce from the hypotheses, Corollary \ref{positivity2} and Proposition \ref{pp2c} that for all $\psi \in {\cal H}$ such that $E_{\Lambda}(U)\psi \in {\cal H}_c(U)\cap {\cal D}(A_{if\nabla\bar{f}}),$
\begin{eqnarray}
c_{\Lambda,f,|\nabla f|^2}^{\sharp} \|E_{\Lambda}(U)\psi\|^2 &\leq & \liminf_{n\rightarrow \pm \infty} n^{-1} \bra U^n E_{\Lambda}(U)\psi, A_{if\nabla\bar{f}} U^n E_{\Lambda}(U) \psi \ket \nonumber \\
&\leq & \limsup_{n\rightarrow \pm \infty} n^{-1} \bra U^n E_{\Lambda}(U)\psi, A_{if\nabla\bar{f}} U^n E_{\Lambda}(U)\psi \ket \leq C_{\Lambda,f,|\nabla f|^2}^{\flat} \|E_{\Lambda}(U)\psi\|^2 \enspace . \label{liminf(a)}
\end{eqnarray}
Now, we follow the proof and notations of Proposition \ref{laurent-2ub}. By Lemma \ref{H0-4} (resp. by hypothesis), $L_f$ (resp. $U$) belongs to $\cap_{j=1}^d C^1(X_j)$. The operators $W:=L_f^*U$ and $U^n$ also belong to $\cap_{j=1}^d C^1(X_j)$ for any $n\in {\mathbb Z}$. So, for $j\in \{1, \ldots, d\}$, we split $C_{j,U}$ as before: $C_{j,U}= C_{j,L_{f}} +L_{f} C_{j,W}L_{f}^*$, where $C_{j,L_{f}}=L_{if\partial_j \bar{f}}$. For any $n\in {\mathbb Z}$ and any $\psi \in {\cal H}$ such that $E_{\Lambda}(U)\psi \in {\cal D}_X$, one also has that $U^nE_{\Lambda}(U) \psi \in {\cal D}_X$ and
\begin{equation*}
\sum_{j=1}^d \Re \bra C_{j,U} U^n E_{\Lambda}(U) \psi, X_j U^n E_{\Lambda}(U)\psi\ket - \bra U^n E_{\Lambda}(U) \psi, A_{if\nabla\bar{f}} U^n E_{\Lambda}(U) \psi \ket = \Re \left( Q_0(U^nE_{\Lambda}(U) \psi, U^nE_{\Lambda}(U) \psi) \right) \enspace ,
\end{equation*}
where $|Q_0(U^nE_{\Lambda}(U) \psi, U^nE_{\Lambda}(U) \psi)|\leq C \|\psi\|^2$ for some $C>0$ according to Lemma \ref{shortcut}.  It follows from (\ref{liminf(a)}) that for any $\psi \in {\cal H}$ such that $E_{\Lambda}(U)\psi \in {\cal H}_c(U)\cap {\cal D}_X$ (${\cal D}_X\subset {\cal D}(A_{if\nabla\bar{f}})$),
\begin{eqnarray*}
c_{\Lambda,f,|\nabla f|^2}^{\sharp} \|E_{\Lambda}(U)\psi\|^2 &\leq & \liminf_{n\rightarrow \pm \infty} n^{-1} \sum_{j=1}^d \Re \bra C_{j,U} U^nE_{\Lambda}(U)\psi, X_j U^nE_{\Lambda}(U)\psi\ket \\
&\leq & \limsup_{n\rightarrow \pm \infty} n^{-1} \sum_{j=1}^d \Re \bra C_{j,U} U^nE_{\Lambda}(U)\psi, X_j U^nE_{\Lambda}(U)\psi\ket \leq C_{\Lambda,f,|\nabla f|^2}^{\flat} \|E_{\Lambda}(U)\psi\|^2 \enspace .
\end{eqnarray*}
We deduce (\ref{liminf(b)}) from Proposition \ref{propag6}. Now, consider $\phi \in C^0({\mathbb T})$ vanishing outside the Borel set $\Lambda\subset {\mathbb T}$. By substituting $c_{\Lambda,f,|\nabla f|^2}^{\sharp}$ to $c_{\Lambda,f,|\nabla f|^2}$, $C_{\Lambda,f,|\nabla f|^2}^{\flat}$ to $C_{\Lambda,f,|\nabla f|^2}$ and $E_{\Lambda}(U)\psi$ to $\Phi(U)\psi$ in the proof of (\ref{liminf(b)}), we can also prove (\ref{liminf(c)}) along the same lines but for any $\psi \in {\cal H}$ such that $\Phi(U)\psi \in {\cal H}_c(U)\cap {\cal D}_X$. The last step consists in using Proposition \ref{10} and noting that $\Phi(U)$ belongs to $\cap_{j=1}^d C^1(X_j)$ whenever $U\in \cap_{j=1}^d C^1(X_j)$, $\phi \in C^1({\mathbb T})$ and $\phi'\in {\cal A}({\mathbb T})$. For any such function $\phi$, we have that for any $n\in {\mathbb Z}$, $\Phi(U)U^n\psi \in {\cal D}_X$ if $\psi \in {\cal D}_X$. This completes the proof of (\ref{liminf(c)}). The last statement follows from the last comment in Corollary \ref{positivity2}.

\subsection{Proof of Proposition \ref{laurent-2tr}}

According to (\ref{mourre-tr}), if $U\in C^1(\sum_{j=1}^d \alpha_j X_j)$ and $\mathrm{ad}_{\sum_{j=1}^d \alpha_j X_j}T^{-\alpha}U$ is compact, then $(U^* (\sum_{j=1}^d \alpha_j X_j) U-(\sum_{j=1}^d \alpha_j X_j)) \simeq |\alpha|^2$. Let $(\varphi,\psi)\in {\cal H}\times {\cal D}(\sum_{j=1}^d \alpha_j X_j)$. If either $\varphi$ or $\psi$ belongs to ${\cal H}_c(U)$, then we get from Lemma \ref{pp2c0},
\begin{equation}\label{lim(a)}
\lim_{n\rightarrow \pm \infty} n^{-1} \bra U^n \psi, (\sum_{j=1}^d \alpha_j X_j) U^n \psi \ket = |\alpha|^2 \bra \varphi,\psi \ket \enspace .
\end{equation}
Combining statement (\ref{lim(a)}) with the strategy developed in Sections 5.4 and 5.5, we obtain that for any $\psi \in {\cal D}_X\cap {\cal H}_c(U),$
\begin{equation*}
\lim_{n\rightarrow \pm \infty} n^{-1} \sum_{j=1}^d \Re \bra C_{j,U} U^n\psi, X_j U^n\psi\ket= |\alpha|^2 \|\psi\|^2 \enspace ,
\end{equation*}
which proves (\ref{lim(b)}) in view of Proposition \ref{propag6}.

\section{On Regularity Classes}\label{tech1}

We have gathered some elementary properties of the regularity classes $C^k(A)$ and ${\cal C}^{s,p}(A)$ that are used throughout the paper. For more details see Chapter 5 in \cite{abmg} or \cite{sah}. In the following, ${\cal H}$ denotes a fixed infinite dimensional Hilbert space, while $A$ is a fixed self-adjoint operator (densely) defined on ${\cal H}$ with domain ${\cal D}(A)$. ${\cal S}$ denotes a dense linear subspace of ${\cal H}$.

\subsection{Main Properties}

In section 2.2, the regularity of a bounded operator defined on ${\cal H}$ w.r.t $A$ is expressed in terms of the derivation defined on ${\cal B}({\cal H})$ by the operation $\mathrm{ad}_A$. From a theoretical point of view, it is often convenient to have in mind the equivalence described by Proposition \ref{equiv0}. Given $B\in {\cal B}({\cal H})$, we define the strongly continuous function ${\cal W}_B$ by:
\begin{eqnarray*}
{\cal W}_B: {\mathbb R} & \rightarrow & {\cal B}({\cal H})\\
t & \mapsto & e^{iAt}Be^{-iAt} \enspace .
\end{eqnarray*}
\begin{prop}\label{equiv0} Let $k\in {\mathbb N}$. The following statements are equivalent:
\begin{itemize}
\item $B\in C^k(A)$
\item The map ${\cal W}_B$ is $C^k$ with respect to the strong topology on ${\cal B}({\cal H})$.
\item The map ${\cal W}_B$ is $C^k$ with respect to the weak topology on ${\cal B}({\cal H})$.
\end{itemize}
Moreover, ${\cal W}_B^{(k)} (0) = i^k \mathrm{ad}_A^k B$.
\end{prop}
See Proposition 5.1.2 and Lemma 6.2.9 in \cite{abmg}. Alternative characterizations are developed in Section 6.2. According to the notations of Proposition \ref{equiv0}, ${\cal B}({\cal H})=C^0(A)$ and for any $B\in {\cal B}({\cal H})$, ${\cal W}_B(0) = B=\mathrm{ad}_A^0 B$. For all nonnegative integers $k$, $C^k(A)$ is clearly a vector subspace of ${\cal B}({\cal H})$ and $C^{k+1}(A) \subset C^k(A)$. We have that:

\begin{prop}\label{2} If $B\in C^1(A)$, then $B({\cal D}(A)) \subset {\cal D}(A)$.
\end{prop}

\begin{prop}\label{3} Let $k\in {\mathbb N}$ and $(B,C) \in C^k(A)\times C^k(A)$. then,
\begin{itemize}
\item $B^* \in C^k (A)$ and for all $j\in \{0,\ldots,k\}$, $\mathrm{ad}_A^j B^* = (-1)^j (\mathrm{ad}_A^j B)^*$
\item $BC \in C^k(A)$ and for all $j\in \{1,\ldots,k\},$
\begin{equation*}
\mathrm{ad}_A^j BC = \sum_{l_1+l_2=j} \frac{j!}{l_1! l_2!} (\mathrm{ad}_A^{l_1} B) (\mathrm{ad}_A^{l_2} C) \enspace.
\end{equation*}
In particular, $\mathrm{ad}_A BC = (\mathrm{ad}_A B) C + B (\mathrm{ad}_A C)$
\item for all $j\in \{0,\ldots, k\}$, $\mathrm{ad}_A^j B \in C^{k-j}(A)$.
\item If $B$ is invertible (i.e $B^{-1}\in {\cal B}({\cal H})$) and $B\in C^1(A)$, then $B^{-1} \in C^1(A)$: $\mathrm{ad}_A B^{-1} = - B^{-1} (\mathrm{ad}_A B) B^{-1}$.
\end{itemize}
\end{prop}
See Propositions 5.1.2, 5.1.5, 5.1.6, 5.1.7 in \cite{abmg} for a proof. Combining the last statements of Proposition \ref{3}, we deduce that if an invertible bounded operator $B$ belongs to $C^k(A)$, then its inverse $B^{-1}$ also belongs to $C^k(A)$.

As anticipated in Section 2, it also possible to deal with intermediate scales of regularity:
\begin{defin}\label{csp} Given a Hilbert space ${\cal H}$, let $B\in {\cal B}({\cal H})$ and $A$ be a self-adjoint operator, densely defined on ${\cal H}$. For $s>0$, $p\in [1,\infty)$, we say that $B$ belongs to the class ${\cal C}^{s,p}(A)$ if there exists $l>s$ such that:
\begin{equation}\label{int}
\int_{-1}^1 \| \sum_{m=0}^l (-1)^m \binom{l}{m}e^{imA\tau}Be^{-imA\tau} \|^p \, \frac{d\tau}{|\tau|^{sp+1}}  < \infty \enspace.
\end{equation}
\end{defin}
The classes ${\cal C}^{s,p}(A)$($s>0$, $p\in [1,\infty)$) are also vector subspaces of ${\cal B}({\cal H})$. In addition, we have that:
\begin{prop}\label{4} Let $(s,p)\in {\mathbb N}$ and $(B,C) \in {\cal C}^{s,p}(A)\times {\cal C}^{s,p}(A)$. then,
\begin{itemize}
\item $B^* \in {\cal C}^{s,p}(A)$
\item $BC \in {\cal C}^{s,p}(A)$
\item If $s=k+\sigma$ with $k\in {\mathbb N}$, $\sigma\in (0,1]$ and $j\in \{0,\ldots, k\}$, $\mathrm{ad}_A^j B \in {\cal C}^{s-j,p}(A)$.
\item If $B$ is invertible, i.e $B^{-1}\in {\cal B}({\cal H})$, then $B^{-1} \in {\cal C}^{s,p}(A)$.
\end{itemize}
\end{prop}
See Propositions 5.2.2, 5.2.3, 5.2.4 in \cite{abmg} for a proof.

The relationships between these regularity classes with the self-adjoint functional calculus are deeply analyzed in \cite{abmg} (Theorem 6.2.5, Corollary 6.2.6) and \cite{gj}. In the unitary context, we also have that:
\begin{prop}\label{10} Let $k\in {\mathbb N}$ and $\phi\in C^0({\mathbb T})$. Assume that $U \in C^k(A)$ and that $(m^k \hat{\phi}_m)_{m\in {\mathbb Z}}\in l^1({\mathbb Z})$ where $(\hat{\phi}_m)_{m\in {\mathbb Z}}$ denotes the sequence of Fourier coefficients of $\phi$. Then, $\Phi(U) \in C^k(A)$ and for all $j\in \{0,\ldots,k\},$
\begin{equation*}
\mathrm{ad}_A^j \Phi(U) = \sum_{m\in {\mathbb Z}} \hat{\phi}_m \mathrm{ad}_A^j U^m \enspace.
\end{equation*}
\end{prop}
{\bf Remark:} If $\phi\in C^k({\mathbb T})$ and $\phi^{(k)}\in {\cal A}({\mathbb T})$, then $(m^k \hat{\phi}_m)_{m\in {\mathbb Z}}\in l^1({\mathbb Z})$ \cite{ed}.

The proof of Proposition \ref{10} is based on the following lemmata whose proofs can also be found in \cite{ggm}:
\begin{lem}\label{H0-1} Let $(B_n)_{n\in {\mathbb N}}\subset C^1(A)$. Assume that the sequences $(B_n)$ and $(\mathrm{ad}_A B_n)$ converge weakly to some (bounded) operators $C_0$ and $C_1$ respectively. Then, $C_0 =w-\lim B_n$ belongs to $C^1(A)$ and
\begin{equation*}
\mathrm{ad}_A C_0 =C_1 \enspace .
\end{equation*}
\end{lem}
\noindent{\bf Proof:} It follows from the hypotheses that the sequences $(B_n)_{n\in {\mathbb N}}$ and $(\mathrm{ad}_A B_n)_{n\in {\mathbb N}}$ are uniformly bounded. In particular, there exists $C>0$ such that for all $(\varphi,\psi)\in {\cal D}(A)\times {\cal D}(A)$ and all $n\in {\mathbb N},$
\begin{equation}\label{lim}
\left| \bra A\varphi,B_n\psi \ket - \bra B_n\varphi,A\psi\ket \right| \leq C \|\varphi\| \|\psi\| \enspace.
\end{equation}
On the other hand, for all $(\varphi,\psi)\in {\cal D}(A)\times {\cal D}(A),$
\begin{eqnarray*}
\bra A\varphi,C_0\psi \ket - \bra C_0\varphi,A\psi\ket &=& \lim_{n\rightarrow \infty} \bra A\varphi,B_n \psi \ket - \bra B_n \varphi,A\psi\ket\\
&=& \lim_{n\rightarrow \infty} \bra \varphi, \mathrm{ad}_A B_n \psi \ket = \bra \varphi, C_1 \psi \ket\enspace.
\end{eqnarray*}
By taking the limit in (\ref{lim}), we deduce that the operator $C_0 \in C^1(A)$ and $\mathrm{ad}_A C_0 = \lim_{n\rightarrow \infty}\mathrm{ad}_A B_n$. \ep

By induction, we get:
\begin{lem}\label{H0-2} Let $k\in {\mathbb N}$, $(B_n)_{n\in {\mathbb N}}\subset C^k(A)$. Assume that for all $j\in \{0,\ldots,k\}$, the sequences $(\mathrm{ad}_A^j B_n)$ converge weakly to some (bounded) operator $C_j$ respectively. Then, $C_0 =w-\lim_{n\rightarrow \infty} B_n$ belongs to $C^k(A)$ and for all $j\in \{0,\ldots,k\},$
\begin{equation*}
\mathrm{ad}_A^j C_0 =C_j \enspace .
\end{equation*}
\end{lem}

\begin{lem} \label{9} Let $k\in {\mathbb N}$. Assume that $U\in C^k(A)$. Then, for all $n\in {\mathbb Z}$, $U^n \in C^k(A)$. Moreover, there exists $C>0$ such that for all $j\in \{0,\ldots,k\}$ and all $|n| \geq j,$
\begin{equation*}
\| \mathrm{ad}_A^j U^n\| \leq C^j |n|^j \enspace,
\end{equation*}
where $C:=\sqrt{\sum_{j=0}^{k} \|\mathrm{ad}_A^j U \|^2}\geq \|\mathrm{ad}_A^0 U \|=\|U \|=1$.
\end{lem}
\noindent{\bf Proof:} The first part is a consequence of Proposition \ref{3}. To prove the second part, it is enough to consider the case of positive $n$ (see again Proposition \ref{3}). We have that:
\begin{equation*}
\mathrm{ad}_A^j U^n = \sum_{(j_1,\ldots,j_n)\in \{0,\ldots,j\}^n, j_1+\ldots+j_n=j} \binom{j}{j_1 \ldots j_n} (\mathrm{ad}_A^{j_1} U)\ldots (\mathrm{ad}_A^{j_n} U) \enspace.
\end{equation*}
Let $n\geq j$. For each term involved in the sum on the RHS, at least $(n-j)$ indices $j_n$'s are zero. Since $\|\mathrm{ad}_A^0 U \|=\|U \|=1$, each one of these terms can be estimated by $C^j$. Since the sum on the RHS involves $n^j$ terms, the estimate follows. \ep
\\

\noindent{\bf Proof of Proposition \ref{10}:} We have that: $\Phi(U) =\sum_{m\in {\mathbb Z}} \hat{\phi}_m U^m$ where the series on the RHS is norm convergent. Since $(m^k \hat{\phi}_m)_{m\in {\mathbb Z}}\in l^1({\mathbb Z})$, it follows from Lemma \ref{9} that for all $j\in \{0,\ldots,k\}$, the series $\sum_m \hat{\phi}_m \mathrm{ad}_A^j U^m$ is norm convergent. The conclusion follows by applying Lemma \ref{H0-2}. \ep

\begin{lem}\label{ad-compact} If $B$ is a compact operator defined on ${\cal H}$ which belongs to ${\cal C}^{1,1}(A)$, then $\mathrm{ad}_A B$ is also compact.
\end{lem}
The proof is the remark (ii) made in the proof of Theorem 7.2.9 in \cite{abmg}. Due to the inclusions (5.2.10) noted in \cite{abmg}, $\mathrm{ad}_A B$ can be expressed as the norm-limit of the family of compact operators $(-i\varepsilon^{-1}(e^{iA\varepsilon}B e^{-iA\varepsilon}-B))_{\varepsilon >0}$ when $\varepsilon$ tends to 0.

\subsection{Some Equivalences}

We conclude this paper by various characterizations of the class $C^1(A)$ for invertible operators. These results build over Lemma 2.1 in \cite{abc1}. We say that $V\in {\cal B}({\cal H})$ is invertible if $V$ is boundedly invertible, i.e. its inverse $V^{-1}$ belongs to ${\cal B}({\cal H})$.

As a general remark, we recall that a bounded invertible operator $V$ belongs to $C^1(A)$ if and only if $V^{-1}$ belongs to $C^1(A)$. Then, for any $n\in {\mathbb Z}$, $V^n\in C^1(A)$ and $V^n {\cal D}(A)= {\cal D}(A)$  (see e.g. Propositions \ref{2} and \ref{3}).

\begin{lem}\label{propag01} Let $V\in C^1(A)$ be a bounded invertible operator and consider the sesquilinear forms $(F_n)_{n\in {\mathbb Z}}$ defined on ${\cal D}(A)\times {\cal D}(A)$ by $F_n(\varphi,\psi):=\langle \varphi,V^{-n}AV^n\psi\rangle\ - \langle\varphi,A\psi\rangle$ if $n\geq 0$ and $F_n(\varphi,\psi):=\langle\varphi,A\psi\rangle-\langle \varphi,V^{-n}AV^n\psi\rangle$ if $n<0$. Then, the forms $(F_n)_{n\in {\mathbb Z}}$ are continuous for the topology induced by ${\cal H}\times {\cal H}$: namely, for all $(\varphi,\psi)\in {\cal D}(A)\times {\cal D}(A)$ and all $n\in {\mathbb Z},$
\begin{equation*}
\left| F_n(\varphi,\psi) \right| \leq C |n| B^{2|n|-1} \|\varphi\| \|\psi\|\enspace ,
\end{equation*}
where $B:=\max \|V^{\pm 1}\|$ and $C:=\max \|\mathrm{ad}_A V^{\pm 1}\|$. Each form $F_n$ extend continuously as a bounded sesquilinear form $F_n^o$ on ${\cal H}\times {\cal H}$. For $n\geq 0$ (resp. $n<0$), denote by $(V^{-n}AV^n-A)^o$ (resp. $(A-V^{-n}AV^n)^o$) the (unique) bounded linear operator associated to $F_n^o$. Then, we have that: $(V^{-n}AV^n-A)^o=V^{-n} (\mathrm{ad}_A V^n)$ (resp. $(A-V^{-n}AV^n)^o=-V^{-n} (\mathrm{ad}_A V^n)$).
\end{lem}
\noindent{\bf Proof:} If $V$ belongs to $C^1(A)$, then for any $n\in {\mathbb Z}$, $(V^n)^*=(V^*)^n$ also belong to $C^1(A)$ and $(V^n)^* {\cal D}(A)= {\cal D}(A)$. The proof for the case $n=0$ is clear. The proofs for $n>0$ and $n<0$ are similar. We restrict our discussion to the case $n>0$. Given any $(\varphi,\psi)\in {\cal D}(A)\times {\cal D}(A),$
\begin{eqnarray*}
F_n(\varphi,\psi) &=& \bra (V^{-n})^*\varphi, AV^n\psi \ket - \bra (V^n)^* (V^{-n})^*\varphi, A \psi \ket\\
                  &=& \bra A(V^{-n})^*\varphi, V^n\psi \ket - \bra (V^{-n})^*\varphi, V^n A \psi \ket = \bra (V^{-n})^*\varphi, (\mathrm{ad}_A V^n)\psi \ket \\
\left| F_n(\varphi,\psi) \right| &\leq & \|\mathrm{ad}_A V^n\| \, B^{|n|} \|\varphi \|\,  \|\psi\| \enspace .
\end{eqnarray*}
Given $n\in {\mathbb N}$, $\mathrm{ad}_A V^n = \sum_{k=0}^{n-1} V^k \left(\mathrm{ad}_A V \right) V^{n-1-k}$. It follows that $\|\mathrm{ad}_A V^n\|\leq C |n| B^{|n|-1}$, hence the first statement of the lemma. We have that for all $(\varphi,\psi)\in {\cal D}(A)\times {\cal D}(A)$, $F_n(\varphi,\psi) = \bra (V^{-n})^* \varphi, (\mathrm{ad}_A V^n)\psi \ket = \bra \varphi, V^{-n}(\mathrm{ad}_A V^n)\psi \ket$. Extending this identity by continuity on ${\cal H}\times {\cal H}$ implies the second part of the Lemma. \ep

\noindent{\bf Remark:} The continuity property introduced in Lemma \ref{propag01} also implies that for any $(\varphi,\psi)\in {\cal H}\times {\cal D}(A)$: $F_n^o(\varphi,\psi) =\langle \varphi,(V^{-n}AV^n-A)^o\psi\rangle =\langle \varphi,V^{-n}AV^n\psi\rangle\ - \langle\varphi,A\psi\rangle$ if $n\geq 0$ and $F_n^o(\varphi,\psi)=\langle \varphi,(A-V^{-n}AV^n)^o\psi\rangle =\langle\varphi,A\psi\rangle-\langle \varphi,V^{-n}AV^n\psi\rangle$ if $n<0$. We deduce that:
\begin{prop}\label{propag07} Let $V\in C^1(A)$ be invertible. Then, for any $(\varphi,\psi) \in {\cal H}\times {\cal D}(A)$ and any $n\in {\mathbb N},$
\begin{eqnarray*}
\langle \varphi,V^{-n}AV^n\psi\rangle\ - \langle\varphi,A\psi\rangle &=& \sum_{m=0}^{n-1} \bra \varphi, V^{-m} (V^{-1}AV-A)^o V^m \psi \ket \\
\langle \varphi,V^n AV^{-n}\psi\rangle\ - \langle\varphi,A\psi\rangle &=& - \sum_{m=1}^{n} \bra \varphi, V^m (V^{-1}AV-A)^o V^{-m} \psi \ket
\end{eqnarray*}
\end{prop}
A similar result can be drawn replacing $(V^{-1}AV-A)^o$ by $(A-VAV^{-1})^o$ in Proposition \ref{propag07}.

Now, we establish various converse statements to Lemma \ref{propag01}.
\begin{lem}\label{propag02} Let $V\in {\cal B}({\cal H})$ be invertible and ${\cal S}$ be a core for $A$. Assume that $V{\cal S}\subset {\cal S}$ and that the sesquilinear form $G_+$ defined on ${\cal S}\times {\cal S}$ by $G_+(\varphi,\psi):=\langle \varphi,V^{-1}AV\psi\rangle\ - \langle\varphi,A\psi\rangle$ is continuous for the topology induced by ${\cal H}\times {\cal H}$. Then, $V$ (and $V^{-1}$) belongs to $C^1(A)$. If we denote by $G_+^o$ the continuous extension of $G_+$ to ${\cal H}\times {\cal H}$ and by $B_+$ the (unique) bounded operator associated to $G_+^o$, we get $VB_+=\mathrm{ad}_A V$.
\end{lem}
\noindent{\bf Proof:} Note first that for any $(\varphi,\psi)\in {\cal H}\times {\cal S}$, $G_+^o(\varphi,\psi) =\langle \varphi,B_+\psi\rangle =\langle \varphi,V^{-1}AV\psi\rangle - \langle\varphi,A\psi\rangle$. We have that for all $(\varphi, \psi)\in {\cal S}\times {\cal S},$
\begin{eqnarray}\label{com1}
\bra A\varphi, V\psi\ket -\bra \varphi, VA\psi\ket &=& \bra (V^{-1})^* V^*\varphi, AV\psi\ket -\bra V^*\varphi, A\psi\ket = \bra V^*\varphi, V^{-1}AV\psi\ket -\bra V^*\varphi, A\psi\ket \nonumber \\
&=& G_+^o(V^*\varphi,\psi)\enspace .
\end{eqnarray}
This implies that: $|\bra A\varphi, V\psi\ket -\bra \varphi, VA\psi\ket | \leq C \|\varphi\| \|\psi\|$ for some $C>0$ and all $(\varphi, \psi)\in {\cal S}\times {\cal S}$. Since ${\cal S}$ is a core for $A$, these estimates extend over ${\cal D}(A)\times {\cal D}(A)$ hence $V\in C^1(A)$. So, (\ref{com1}) rewrites: $\bra \varphi, (\mathrm{ad}_A V)\psi\ket = \bra V^*\varphi, B_+ \psi\ket  = \bra \varphi, VB_+ \psi\ket$ for all $(\varphi, \psi)\in {\cal D}(A)\times {\cal D}(A)$. Extending this identity by continuity on ${\cal H}\times {\cal H}$ concludes the proof. \ep 

\noindent{\bf Remark:} Under the hypotheses of Lemma \ref{propag02}, we also deduce that $\mathrm{ad}_A V^{-1}=-V^{-1}(\mathrm{ad}_A V) V^{-1}= -B_+ V^{-1}$ (see e.g. Proposition \ref{3}).

If ${\cal S}={\cal D}(A)$ in Lemma \ref{propag02}, the form $G_+$ coincides with the form $F_1$ defined in Lemma \ref{propag01}. In this case, we would recover the identity $V(V^{-1}AV-A)^o=\mathrm{ad}_A V$. Actually, we have that:

\begin{lem}\label{propag05} Let $V\in {\cal B}({\cal H})$ be invertible and ${\cal S}$ be a core for $A$. Assume that $V{\cal S}\subset {\cal S}$ and that the sesquilinear form $G_+$ defined on ${\cal S}\times {\cal S}$ by $G_+(\varphi,\psi):=\langle \varphi,V^{-1}AV\psi\rangle\ - \langle\varphi,A\psi\rangle$ is continuous for the topology induced by ${\cal H}\times {\cal H}$. Let us denote by $G_+^o$ the continuous extension of $G_+$ to ${\cal H}\times {\cal H}$ and by $B_+$ the bounded operator associated to $G_+^o$. Then,
\begin{itemize}
\item $V$ (and $V^{-1}$) belongs to $C^1(A)$,
\item $V{\cal D}(A) \subset {\cal D}(A)$ and the sesquilinear form $F_1$ defined on ${\cal D}(A)\times {\cal D}(A)$ by $F_1(\varphi,\psi):=\langle \varphi,V^{-1}AV\psi\rangle\ - \langle\varphi,A\psi\rangle$ is continuous for the topology induced by ${\cal H}\times {\cal H}$,
\item if $F_1^o$ and $(V^{-1}AV-A)^o$ denote respectively the continuous extension of $F_1$ to ${\cal H}\times {\cal H}$ and the bounded operator associated to $F_1^o$, we get: $F_1^o=G_+^o$ and $B_+=(V^{-1}AV-A)^o$.
\end{itemize}
In particular, $VB_+=V(V^{-1}AV-A)^o=\mathrm{ad}_A V$.
\end{lem}
\noindent{\bf Proof:} The first statement quotes Lemma \ref{propag02}. The remaining parts follow from Lemma \ref{propag01} once observed that $G_+=F_1\upharpoonright {\cal S}\times {\cal S}$ and that ${\cal S}\times {\cal S}$ is dense in ${\cal H}\times {\cal H}$. \ep

Permuting the roles of $V$ and $V^{-1}$, we also have that:
\begin{lem}\label{propag02bis} Let $V\in {\cal B}({\cal H})$ be invertible and ${\cal S}$ be a core for $A$. Assume that $V^{-1}{\cal S}\subset {\cal S}$ and that the sesquilinear form $G_-$ defined on ${\cal S}\times {\cal S}$ by $G_-(\varphi,\psi):= \langle\varphi,A\psi\rangle -\langle \varphi,VAV^{-1}\psi\rangle$ is continuous for the topology induced by ${\cal H}\times {\cal H}$. Then, $V^{-1}$ (and $V$) belongs to $C^1(A)$. If we denote by $G_-^o$ the continuous extension of $G_-$ to ${\cal H}\times {\cal H}$ and by $B_-$ the (unique) bounded operator associated $G_-^o$, we get $V^{-1}B_-=-\mathrm{ad}_A V^{-1}$.
\end{lem}
\noindent{\bf Remark:} Under the hypotheses of Lemma \ref{propag02bis}, we also deduce that $\mathrm{ad}_A V=-V(\mathrm{ad}_A V^{-1}) V= B_-V$.

If ${\cal S}={\cal D}(A)$ in Lemma \ref{propag02bis}, the form $G_-$ coincides with $F_{-1}$ defined in Lemma \ref{propag01}. In this case, we would recover the identity $V^{-1}(A-VAV^{-1})^o=-\mathrm{ad}_A V^{-1}$. Actually, we have that:

\begin{lem}\label{propag05bis} Let $V\in {\cal B}({\cal H})$ be invertible operator and ${\cal S}$ be a core for $A$. Assume that $V^{-1}{\cal S}\subset {\cal S}$ and that the sesquilinear form $G_-$ defined on ${\cal S}\times {\cal S}$ by $G_-(\varphi,\psi):= \langle\varphi,A\psi\rangle -\langle \varphi,VAV^{-1}\psi\rangle$ is continuous for the topology induced by ${\cal H}\times {\cal H}$. Let us denote by $G_-^o$ the continuous extension of $G_-$ to ${\cal H}\times {\cal H}$ and by $B_-$ the bounded operator associated to $G_-^o$. Then,
\begin{itemize}
\item $V^{-1}$ (and $V$) belongs to $C^1(A)$,
\item $V^{-1}{\cal D}(A) \subset {\cal D}(A)$ and the sesquilinear form $F_{-1}$ defined on ${\cal D}(A)\times {\cal D}(A)$ by $F_{-1}(\varphi,\psi):=\langle\varphi,A\psi\rangle -\langle \varphi,VAV^{-1}\psi\rangle$ is continuous for the topology induced by ${\cal H}\times {\cal H}$,
\item if $F_{-1}^o$ and $(A-VAV^{-1})^o$ denote respectively the continuous extension of $F_{-1}$ to ${\cal H}\times {\cal H}$ and the bounded operator associated to $F_{-1}^o$, we get: $F_{-1}^o=G_-^o$ and $B_-=(A-VAV^{-1})^o$.
\end{itemize}
In particular, $V^{-1}(A-VAV^{-1})^o=V^{-1}B_- =-\mathrm{ad}_A V^{-1}$.
\end{lem}

\noindent{\bf Remark:} Assuming that $V\in C^1(A)$ is invertible, we have that $V^{-1}V=I=VV^{-1}$, $(\mathrm{ad}_A V) V^{-1}+ V(\mathrm{ad}_A V^{-1})=0=(\mathrm{ad}_A V^{-1}) V+ V^{-1}(\mathrm{ad}_A V)$. We deduce from Lemma \ref{propag01} that: $V^{-1}(A-VAV^{-1})^oV=(V^{-1}AV-A)^o$.

The next characterization of the class $C^1(A)$ is proved in Proposition 2.2 \cite{ggm}:
\begin{lem}\label{propag03} An operator $B\in {\cal B}({\cal H})$ belongs to $C^1(A)$ if and only if $B {\cal D}(A)\subset {\cal D}(A)$ and the operator $AB-BA:{\cal D}(A)\rightarrow {\cal H}$ extends to a bounded operator on ${\cal H}$. In this case, $AB=BA+\mathrm{ad}_A B$ as an identity on ${\cal D}(A)$.
\end{lem}

For invertible operators, we also have that:
\begin{lem}\label{propag04} If $V\in C^1(A)$ is invertible, then $V {\cal D}(A)\subset {\cal D}(A)$ and the operator $V^{-1}AV-A:{\cal D}(A)\rightarrow {\cal H}$ extends (uniquely) to a bounded operator on ${\cal H}$. In this case, $V^{-1}AV=A+(V^{-1}AV-A)^o= A+V^{-1}\mathrm{ad}_A V=A-(\mathrm{ad}_A V^{-1})V$ as an identity on ${\cal D}(A)$.
\end{lem}
\noindent{\bf Proof:} The first statement quotes Proposition \ref{2}. We deduce from the remark following Lemma \ref{propag01} that: for all $(\varphi, \psi)\in {\cal H}\times {\cal D}(A)$, $\bra \varphi, V^{-1}AV\psi\ket - \bra \varphi, A\psi\ket =\bra \varphi, (V^{-1}AV-A)^o\psi\ket$ where $(V^{-1}AV-A)^o$ is bounded. This implies our claim. \ep

\begin{lem}\label{propag06} Let ${\cal S}$ be a core for $A$ such that $V{\cal S} \subset {\cal S}$. Assume that the operator $V^{-1}AV-A:{\cal S}\rightarrow {\cal H}$ extends to a bounded operator on ${\cal H}$ and denote by $C_+$ this extension. Then, $V$ (and $V^{-1}$) belongs to $C^1(A)$ and $\mathrm{ad}_A V=VC_+$.
\end{lem}
\noindent{\bf Proof:} For all $(\varphi, \psi)\in {\cal S}\times {\cal S}$, $\bra A\varphi, V\psi\ket -\bra \varphi, VA\psi\ket = \bra (V^{-1})^*V^*\varphi, AV\psi\ket -\bra V^*\varphi, A\psi\ket = \bra V^*\varphi, V^{-1}AV\psi\ket -\bra V^*\varphi, A\psi\ket = \bra V^*\varphi, C_+\psi\ket = \bra \varphi, V C_+\psi\ket$. The identity extends continuously over ${\cal D}(A) \times {\cal D}(A)$. This clearly shows that $V\in C^1(A)$ and that: $\mathrm{ad}_A V=VC_+$. \ep

Since a bounded invertible operator and its inverse belong simultaneously to the class $C^1(A)$, we have that:
\begin{lem}\label{propag04bis} If $V\in C^1(A)$ is invertible, then $V^{-1}{\cal D}(A)\subset {\cal D}(A)$ and the operator $A-VAV^{-1}:{\cal D}(A)\rightarrow {\cal H}$ extends (uniquely) to a bounded operator on ${\cal H}$. In this case, $VAV^{-1}=A- (A-VAV^{-1})^o= A+V \mathrm{ad}_A V^{-1}= A- (\mathrm{ad}_A V)V^{-1}$ as an identity on ${\cal D}(A)$.
\end{lem}

\begin{lem}\label{propag06bis} Let ${\cal S}$ be a core for $A$ such that $V^{-1}{\cal S} \subset {\cal S}$. Assume that the operator $A-VAV^{-1}:{\cal S}\rightarrow {\cal H}$ extends to a bounded operator on ${\cal H}$ and denote by $C_-$ this extension. Then, $V^{-1}$ (and $V$) belongs to $C^1(A)$ and $\mathrm{ad}_A V^{-1}=-V^{-1}C_-$.
\end{lem}

\noindent{\bf Acknowledgment: } The authors thank J. Asch, A. Joye and the referees for their informative comments.  



\begin{thebibliography}{xxxxxxx}

\bibitem{abmg} W.O. Amrein, A. Boutet de Monvel, V. Georgescu, {\it $C_0$-groups, Commutator Methods and Spectral Theory of Hamiltonians}, Birkh\"auser, 1996.
\bibitem{ak} J. Asch, A. Knauf, {\it Quantum transport on KAM tori}, Comm. Math. Phys. 205 (1999), no. 1, 113-128.
\bibitem{abc1} M.A Astaburuaga, O. Bourget, V.H. Cort\'es, {\it Commutation Relations for Unitary Operators I}, J. Func. Anal. (2015), http://dx.doi.org/10.1016/j.jfa.2015.01.011 
\bibitem{sah} A. Boutet de Monvel, J. Sahbani, {\it On the spectral properties of discrete Schr\"odinger operators: the multi-dimensional case}, Rev. Math. Phys. 11 (1999), no. 9, 1061-1078.
\bibitem{D2} E.B. Davies, {\it Spectral theory and differential operators}, Cambridge, 1995.
\bibitem{ed} R.E. Edwards, {\it Fourier Series: A Modern Introduction}, Vol. I, Springer, 1979.
\bibitem{ev} V. Enss, K. Veselic, {\it Bound states and propagating states for time-dependent Hamiltonians}, Ann. Inst. Henri Poincar\'e, 39 (1983) no. 2, 159-191.
\bibitem{ggm} V. Georgescu, C. G\'erard, J.S. Moller, {\it Commutators, C0-semigroups and resolvent estimates}, J. Funct. Anal. 216 (2004), no. 2, 303-361.
\bibitem{g} J. Geronimus, {\it On the character of the solution of the moment-problem in the case of the periodic in the limit associated fraction} (Russian), Bull. Acad. Sci. URSS. S\'er. Math. [Izvestia Akad. Nauk SSSR] 5 (1941), 203-210.
\bibitem{gt} J.S. Geronimo, A. Teplyaev, {\it A difference equation arising from the trigonometric moment problem having random reflection coefficients - an operator-theoretic approach}, J. Funct. Anal. 123 (1994), no. 1, 12-45.
\bibitem{gj} S. Gol\'enia, T.  Jecko, {\it A new look at Mourre's commutator theory}, Complex Anal. Oper. Theory 1 (2007), no. 3, 399-422.
\bibitem{gnva1} L. Golinskii, P. Nevai, W. Van Assche, {\it Perturbation of orthogonal polynomials on an arc of the unit circle}, J. Approx. Theory 83 (1995), no. 3, 392-422.
\bibitem{gnva2} L. Golinskii, P. Nevai, F. Pint\'er, W. Van Assche, {\it Perturbation of orthogonal polynomials on an arc of the unit circle. II}, J. Approx. Theory 96 (1999), no. 1, 1-32.
\bibitem{gol} L. Golinskii, {\it Operator theoretic approach to orthogonal polynomials on an arc of the unit circle}, Mat. Fiz. Anal. Geom. 7 (2000), no. 1, 3-34.
\bibitem{gn} L. Golinskii, P. Nevai, {\it Szeg\"o difference equations, transfer matrices and orthogonal polynomials on the unit circle}, Comm. Math. Phys. 223 (2001), no. 2, 223-259.
\bibitem{N2} M.G. Nadkarni, {\it Spectral Theory of Dynamical Systems}. Birkh\"auser, 1998.
\bibitem{rs} M. Reed, B. Simon, {\it Methods of Modern Mathematical Physics}, Vol.1-4, Academic Press, 1975-1979.
\bibitem{sah1} J. Sahbani, {\it Spectral theory of certain unbounded Jacobi matrices}, J. Math. Anal. Appl. 342 (2008), no. 1, 663-681.
\bibitem{sah2} J. Sahbani, {\it Spectral properties of Jacobi matrices of certain birth and death processes}, J. Operator Theory 56 (2006), no. 2, 377-390.
\bibitem{sah3} J. Sahbani, {\it Mourre's theory for some unbounded Jacobi matrices}, J. Approx. Theory 135 (2005), no. 2, 233-244.
\bibitem{S1} B. Simon, {\it Orthogonal polynomials on the unit circle. Parts 1 \& 2}, American Mathematical Society, Colloquium Publications 54, 2005.
\bibitem{S2} B. Simon, {\it Szeg\"o's theorem and its descendants. Spectral theory for $L^2$ perturbations of orthogonal polynomials}, M. B. Porter Lectures. Princeton University Press, 2011.
\bibitem{katz} Y. Katznelson, {\it An introduction to harmonic analysis}, Cambridge, 2004.
\bibitem{schmu} K. Schm\"udgen, {\it On domains of powers of closed symmetric operators}, J. Operator Theory 9 (1983), no. 1, 53-75.


\end{thebibliography}
\end{document}